\documentclass[12pt]{article}
\usepackage{amsmath}
\usepackage{amssymb}
\usepackage{latexsym}
\usepackage{amsthm}
\usepackage{mathrsfs}
\usepackage{fancyhdr}
\usepackage{amsfonts}
\usepackage{amssymb}
\usepackage{graphicx}

\newtheorem{thm}{Theorem}[section]
\newtheorem{lemma}[thm]{Lemma}
\newtheorem{cor}[thm]{Corollary}

\newtheorem{prob}[thm]{Open problem}
\newcommand{\pf}{\noindent{Proof:\, }}
\newcommand{\vs}{\vspace{3mm}}

\usepackage{latexsym , bm}

\begin{document}

\title{On $k$-resonant fullerene graphs\footnote{This
 work is  supported by  NSFC (Grant no. 10831001).}}

\author{Dong Ye, Zhongbin Qi and Heping Zhang\footnote{Corresponding author.}}
 \date{\small{School of Mathematics and Statistics, Lanzhou University,
 Lanzhou, Gansu 730000, P. R. China}\\
 \small{E-mails: dye@lzu.edu.cn, qizhb02@st.lzu.edu.cn, zhanghp@lzu.edu.cn}}
 \maketitle


\begin{abstract}
A fullerene graph $F$ is a 3-connected plane cubic graph with
exactly 12 pentagons and the remaining hexagons. Let $M$ be a
perfect matching of $F$. A cycle $C$ of $F$ is $M$-alternating if
the edges of $C$ appear alternately in and off $M$. A set $\mathcal
H$ of disjoint hexagons of $F$ is called a resonant pattern (or
sextet pattern) if $F$ has a perfect matching $M$ such that all
hexagons in $\mathcal H$ are $M$-alternating. A fullerene graph $F$
is $k$-resonant if any $i$ ($0\leq i \leq k$)  disjoint hexagons of
$F$ form a resonant pattern. In this paper, we prove that every
hexagon of a fullerene graph is resonant and all leapfrog fullerene
graphs are 2-resonant. Further, we show that a 3-resonant fullerene
graph has at most 60 vertices and construct all nine 3-resonant
fullerene graphs,  which are also $k$-resonant for every integer
$k>3$. Finally, sextet polynomials of the 3-resonant fullerene
graphs are computed.

\vspace{0.3cm}

\noindent {\em Keywords:} Fullerene graph; Perfect matching;
Resonant pattern; $k$-resonance; Sextet polynomial

\noindent{\em AMS 2000 subject classification:} 05C70, 05C90
\end{abstract}


\section{Introduction}

A {\em fullerene graph} is a  3-connected plane cubic graph with
exactly 12 pentagonal faces and the other faces being hexagonal.
Fullerene graphs have been studied in mathematics as trivalent
polyhedra for a long time \cite{GM,GB},   for example, the
dodecahedron is the fullerene graph with 20 vertices. Fullerene
graphs have been studied in chemistry as fullerene molecules which
have extensive applications in physics, chemistry and material
science \cite{FM}.

Let $G$ be a plane 2-connected graph. A {\em perfect matching} or
{\em 1-factor} $M$ of $G$ is a set of independent edges such that
every vertex of $G$ is incident with exactly one edge in $M$. A
cycle $C$ of $G$ is {\em $M$-alternating} if the edges of $C$ appear
alternately in and off $M$. For a fullerene graph $F$,  every edge
of $F$ belongs to a perfect matching of $F$ \cite{KL,D0}. A hexagon
$h$ of a fullerene graph $F$ is {\em resonant} if $F$ has a perfect
matching $M$ such that $h$ is $M$-alternating. It was proved that
every hexagon of a normal benzenoid system is resonant \cite{ZC}.
This result was generalized to normal coronoid systems \cite{ZZ} and
plane elementary bipartite graphs \cite{ZZ2}. However a fullerene
graph is  a non-bipartite graph. It is natural to ask if every
hexagon of a fullerene graph is resonant. The present paper first
uses Tutte's 1-factor theorem to give a positive answer to this
question.

A set $\mathcal H$ of disjoint hexagons of a fullerene graph $F$ is
a {\em resonant pattern} (or sextet pattern), in other words, such
hexagons are  {\em mutually resonant}, if $F$ has a perfect matching
$M$ such that every hexagon in $\mathcal H$ is $M$-alternating;
equivalently, if $F-\mathcal H$ has a perfect matching, where
$F-\mathcal H$ denotes the subgraph obtained from $F$ by deleting
all vertices of $\mathcal H$ together with their incident edges. The
maximum cardinality of resonant patterns of $F$ is called the {\em
Clar number} of $F$ \cite{Clar}, and the maximum number of
$M$-alternating hexagons over all perfect matchings $M$ of $F$ is
called the {\em Fries number} of $F$ \cite{Fries}. Graver \cite{GJ}
explored some connections among the Clar number, the face
independence number and the Fries number of a fullerene graph, and
obtained a lower bound for the Clar number of leapfrog fullerene
graphs with icosahedral symmetry. Zhang and Ye \cite{ZY1} showed
that the Clar number of a fullerene graph $F_n$ with $n$ vertices
satisfies $c(F_n)\le \lfloor \frac{n-12} 6\rfloor$, which is sharp
for infinitely many fullerene graphs, including $\text{C}_{60}$
whose Clar number is 8 \cite{EB}. Shiu, Lam and Zhang \cite{SLZ1}
computed the Clar polynomial and the sextet polynomial of
$\text{C}_{60}$ by showing that every hexagonal face independent set
of $\text{C}_{60}$ is also a resonant pattern.

A fullerene graph is {\em $k$-resonant} if any $i$ ($0\le i\le k$)
disjoint hexagons are mutually resonant. So $k$-resonant fullerene
graphs are also $(k-1)$-resonant for integer $k\ge 1$. Hence a
fullerene graph with each hexagon being resonant is 1-resonant.
Zheng \cite{Z1, Z2} characterized general $k$-resonant benzenoid
systems. In particular, he showed that every $3$-resonant benzenoid
system is also $k$-resonant ($k\ge 3$). This result also holds for
coronoid systems \cite{CG,LC}, open-ended nanotubes \cite{ZW},
toroidal polyhexes \cite{SLZ2, ZY2} and Klein-bottle polyhexes
\cite{SZ}. For a recent survey on $k$-resonant benzenoid systems,
refer to \cite{G1}.

Here we consider $k$-resonant fullerene graphs. We show that all
leapfrog fullerene graphs are 2-resonant and a 3-resonant fullerene
graph has at most 60 vertices. We construct all 3-resonant fullerene
graphs,  and show that they are all $k$-resonant for every integer
$k\ge 3$. This result is consistent with the aforementioned results.
Finally, sextet polynomials of the 3-resonant fullerene graphs are
computed.
\section{1-resonance of fullerene graphs}
Let $G$ be a plane graph admitting a perfect matching with
vertex-set $V(G)$ and edge-set $E(G)$. Use $\partial G$ denote the
{\em boundary} of $G$, i.e. the boundary of the infinite face of
$G$. For a face $f$ of $G$, let $V(f)$ and $E(f)$ be the sets of
vertices and edges of $f$, respectively. If $G$ is a 2-connected
plane graph, then each face of $G$ is bounded by a cycle. For
convenience, a face is often represented by its boundary if
unconfused. In particular, for a fullerene graph $F$, any {\em
pentagon}, a cycle with length five, and any {\em hexagon}, a cycle
with length six, of $F$ must bound a face since $F$ is cyclically
5-edge connected \cite{D,ZY1}. For a plane graph $G$, a face $f$ of
$G$ {\em adjoins} a subgraph $G'$ of $G$ if $f$ is not a face of
$G'$ and $f$ has an edge in common with $G'$. The faces adjoining
$G'$ are always called {\em adjacent faces} of $G'$. A subgraph $H$
of $G$ is called {\em nice} in \cite{LP} or {\em central} in
\cite{RST} if $G-V(H)$ has a perfect matching. So a resonant pattern
of $G$ can be viewed as a central subgraph of $G$. A graph $G$ is
{\em cyclically $k$-edge connected} if deleting fewer than $k$ edges
of $G$ can not separate $G$ into two components each of which
contains a cycle. By Tutte's Theorem on perfect matchings of graphs
(\cite{LP}, Theorem 3.1.1), we have the following result.

\begin{lemma}\label{thm2-1}
A subgraph $H$ of a graph $G$ is central if and only if for any
$S\subseteq V(G-H)$,
$$C_o(G-H-S)\leq |S|,$$
where $C_o(G-H-S)$ is the number of odd components of $G-H-S$.
\end{lemma}

\begin{thm}\label{thm2-2}
Let $G$ be a cyclically 4-edge connected cubic graph with a 6-length
cycle. Then for every 6-length cycle $H$ of $G$, either $H$ is
central or $G-H$ is bipartite.
\end{thm}

\pf Let $H$ be a 6-length cycle in $G$. If $G-H$ has a perfect
matching, then the theorem holds. If not, then by Lemma \ref{thm2-1}
there exists an $S\subset V(G-H)$ such that $C_o(G-H-S)\ge |S|+2$ by
parity, i.e. $|S|\le C_o(G-H-S)-2$. Since $G$ is cubic, $S$ sends
out at most $3|S|\le 3C_o(G-H-S)-6$ edges.

Let $G_1,G_2,...,G_k$ be all odd components of $G-H-S$, where
$k=C_o(G-H-S)$. Because $G$ is cyclically 4-edge connected and
cubic, it has no cut edges. Every $G_i$ ($i=1,2,...,k$) sends odd
number edges, hence at least three edges, to $H\cup S$. So
$\cup_{i=1}^kG_i$ sends out at least $3C_o(G-H-S)$ edges to either
$S$ or $H$. Since $H$ is a 6-length cycle, there are at most 6 edges
between  $H$ and $\cup_{i=1}^kG_i$. So $\cup_{i=1}^kG_i$ sends at
least $3C_o(G-H-S)-6$ edges to $S$. Hence there are precisely
$3C_o(G-H-S)-6$ edges between $S$ and $\cup_{i=1}^kG_i$. So $S$ is
an independent set, and every $G_i$ sends out exactly 3 edges, and
$G-H-S$ has no even component. In addition, since $G$ is cyclically
4-edge connected, every $G_i$ is a tree. We claim that each $G_i$ is
a singular vertex. If not, then an odd component $G_i$ has at least
2 vertices. So $G_i$ has at least two leaves. Every leaf of $G_i$ is
adjacent to at least two vertices in $S\cup H$. So $G_i$ sends at
least four edges out, contradicting the fact that every $G_i$ sends
precisely three edges out. Therefore $G-H$ is a bipartite graph with
bipartition $(S, V(G-H-S))$. This completes the proof of the
theorem. \qed

\begin{lemma}\cite{D,QZ}\label{lem2-3}
Every fullerene graph is cyclically 5-edge connected.
\end{lemma}

By Lemma \ref{lem2-3} and Theorem \ref{thm2-2}, we immediately have
the following result.

\begin{thm} \label{thm2-4}
Every hexagon of a fullerene graph is resonant.
\end{thm}

\pf Let $F$ be a fullerene graph and $H$ be a hexagon of $F$. It is
obvious that $F-H$ is not bipartite. By Theorem \ref{thm2-2} and
Lemma \ref{lem2-3}, $H$ is central. That means $H$ is resonant. \qed


\section{2-resonant fullerene graphs}

\begin{figure}[!hbtp]\refstepcounter{figure}\label{fig3-1}
\begin{center}
\includegraphics{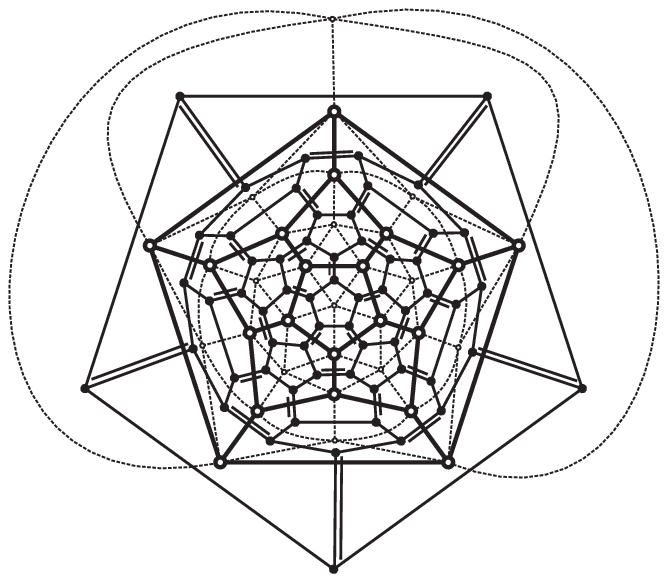}\\
{Figure \ref{fig3-1}: The leapfrog operation on the dodecahedron
$F_{20}$ and the perfect matching $M^0$ of $\text{C}_{60}$ (double
edges). }
\end{center}
\end{figure}

Let $F$ be a fullerene graph. The {\em leapfrog operation} on $F$ is
defined \cite{FP} as follows: for any face $f$ of $F$, add a new
vertex $v_f$ in $f$ and join $v_f$ to all vertices in $V(f)$ to
obtain a new triangular graph $F'$; then take the geometry dual of
the graph $F'$ and denote it by $F^*$ (see Figure \ref{fig3-1}).
Clearly, $F^*$ is a fullerene graph since every vertex of $F'$ is
6-degree excluding exactly 12 5-degree vertices and every face of
$F'$ is a triangle. The edges of $F^*$ cross the edges of $F\subset
F'$ in the geometry dual operation form a perfect matching $M^0$ of
$F^*$. A fullerene graph is called {\em leapfrog fullerene} if it
arises from a fullerene graph by the leapfrog operation. Several
characterizations of leapfrog fullerenes have been given; see Liu,
Klein and Schmalz \cite{LKS}, Fowler and Pisanski \cite{FP}, and
Graver \cite{GJ,GJ2}. For example, a fullerene graph is a leapfrog
fullerene if and only if it has a perfect Clar structure (i.e. a set
of disjoint faces including all vertices); and if and only if it has
a Fries structure (i.e. a perfect matching which avoids edges in
pentagons and is alternating on the maximal number $n/3$ of
hexagons).

Let $F^*$ be a leapfrog fullerene graph arising  from $F$. A face
$f$ of  $F^*$ is called a {\em heritable face} if it lies completely
in some face of $F$, and a {\em fresh face}, otherwise. For example,
$\text{C}_{60}$ is the leapfrog fullerene graph of the dodecahedron
and every pentagon is a heritable face and all hexagons are fresh
faces. The perfect matching $M^0$ corresponds to the Fries structure
of $\text{C}_{60}$ (see Figure \ref{fig3-1}). For a leapfrog
fullerene graph, we have the following result.

\begin{lemma}\label{lem3-1}
Let $F$ be a leapfrog fullerene graph. Then every fresh face is
$M^0$-alternating and all heritable faces are independent.
\end{lemma}

Let $F$ be a leapfrog fullerene and $f$  a heritable face of $F$. A
subgraph of $F$ consisting of $f$ together with all adjacent (fresh)
faces is called the {\em  territory} of $f$, and denoted by $T[f]$.
For two heritable faces $f_1$ and $f_2$, it is easily seen that
there are at most 2 common fresh faces in their territories, which
are adjacent.

\begin{thm}\label{thm3-2}
Every leapfrog fullerene graph is 2-resonant.
\end{thm}

\pf Let $F$ be a leapfrog fullerene graph and $f_1, f_2$ any two
disjoint hexagons. If both $f_1$ and $f_2$ are fresh faces, then
clearly $M^0$ is alternating on both of them by Lemma \ref{lem3-1}.
So suppose that at least one of them is a heritable face, say $f_1$.
Let us denote the six fresh hexagons in $T[f_1]$ by $h_0, h_1
,\dots, h_5$ in clockwise order. If $f_2$ is fresh, then
$f_2\nsubseteq T[f_1]$ and it adjoins at most one of $h_0,h_1,\dots,
h_5$ since $F$ is a leapfrog fullerene graph. If $f_2$ adjoins none
of $h_1,h_3$ and $h_5$, let $M_1:=M^0\oplus h_1\oplus h_3\oplus
h_5$; otherwise, let $M_1:=M^0\oplus h_0\oplus h_2\oplus h_4$. Then
$M_1$ is a perfect matching and alternating on both $f_1$ and $f_2$.
So, in the following, we suppose both $f_1$ and $f_2$ are heritable.
Let $h'_0,h'_1,\dots,h'_5$ be the six fresh hexagons of $T[f_2]$ in
clockwise order. If $T[f_1]$ and $T[f_2]$ have a common hexagon,
then they have exactly two common adjacent hexagons. Assume
$h_{i_0}=h'_{j_0}$ for some $i_0,j_0\in \mathbb Z_6$. Let
$M_2:=M^0\oplus h_{i_0}\oplus h_{i_0+2} \oplus h_{i_0+4}\oplus
h'_{j_0+2} \oplus h'_{j_0+4}$. It is clear that $M_2$ is a perfect
matching alternating on both $f_1$ and $f_2$. Now suppose $T[f_1]$
and $T[f_2]$ have no common hexagons. If no face in $T[f_2]$ adjoins
one of $h_1, h_3 $ and $ h_5$, let $M_3:=M^0\oplus h_1\oplus
h_3\oplus h_5\oplus h'_1\oplus h'_3\oplus h'_5$; otherwise, let
$M_3:=M^0\oplus h_0\oplus h_2\oplus h_4\oplus h_1'\oplus h_3'\oplus
h_5'$. Then $M_3$ is also a perfect matching alternating on both
$f_1$ and $f_2$. So the theorem holds. \qed

\begin{figure}[!hbtp]\refstepcounter{figure}\label{fig3-2}
\begin{center}
\includegraphics{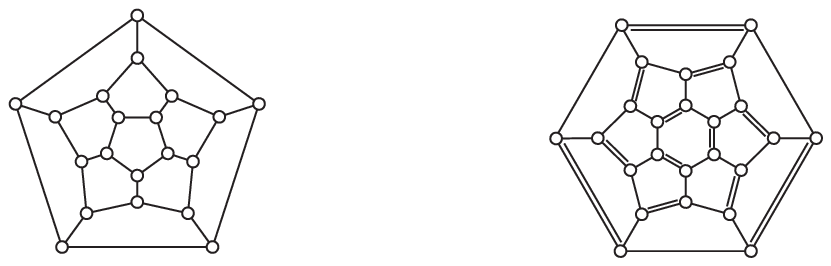}\\
{Figure \ref{fig3-2}: The dodecahedron $F_{20}$ (left) and the
fullerene graph $F_{24}$ with a perfect matching $M$ (right).}
\end{center}
\end{figure}

There exist  2-resonant fullerene graphs which are non-leapfrog. The
dodecahedron $F_{20}$ is a trivial example. The fullerene graph
$F_{24}$, as shown in Figure \ref{fig3-2} (right), is 2-resonant
since the two hexagons  are simultaneously $M$-alternating. Another
non-trivial example is $\text{C}_{70}$.

\begin{lemma}\label{lem3-3}
$\text{C}_{70}$ is 2-resonant.
\end{lemma}

\begin{figure}[!hbtp]\refstepcounter{figure}\label{fig3-3}
\begin{center}
\includegraphics{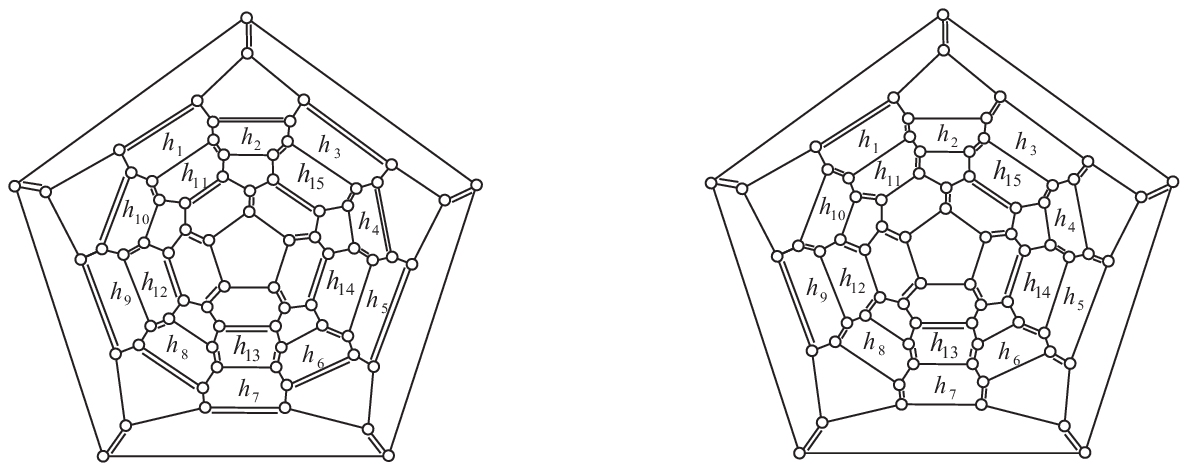}\\
{Figure \ref{fig3-3}: $\text{C}_{70}$ with two perfect matchings
$M_1$ (left) and $M_2$ (right).}
\end{center}
\end{figure}

\pf $\text{C}_{70}$  has two perfect matchings $M_1$ and $M_2$ as
shown in Figure \ref{fig3-3}. It has a total 25 of hexagons. The
hexagons other than $h_1,h_3,h_5,h_7$ and $h_9$ are all
$M_1$-alternating. Let $M_3:=M_1\oplus h_2\oplus h_4\oplus h_6\oplus
h_8\oplus h_{10}$. Then the hexagons other than
$h_{11},h_{12},h_{13},h_{14}$ and $h_{15}$ are all
$M_3$-alternating. We choose any pair of disjoint hexagons $h$ and
$h'$ in $\text{C}_{70}$. If $h,h'\notin \{h_1, h_3, h_5, h_7,h_9\}$,
then $h$ and $h'$ are simultaneously $M_1$-alternating. If
$h,h'\notin \{h_{11},h_{12},h_{13},h_{14},h_{15}\}$, then $h$ and
$h'$ are simultaneously $M_3$-alternating. So suppose $h\in
\{h_1,h_3,h_5,h_7,h_9\}$ and
$h'\in\{h_{11},h_{12},h_{13},h_{14},h_{15}\}$. By symmetry, we may
assume $h=h_1$. If $h'\in \{h_{12}, h_{15}\}$, we may let
$h'=h_{12}$ by the symmetry of $h_{12}$ and $h_{15}$. Then both $h$
and $h'$ are $M_2$-alternating. Finally, if $h'\in \{h_{13},
h_{14}\}$, then $h$ and $h'$ are simultaneously $M_4$-alternating,
where  $M_4:=M_1\oplus h_2\oplus h_{10}$. Hence $\text{C}_{70}$ is
2-resonant. \qed

\begin{figure}[!hbtp]\refstepcounter{figure}\label{fig3-4}
\begin{center}
\scalebox{0.8}{\includegraphics{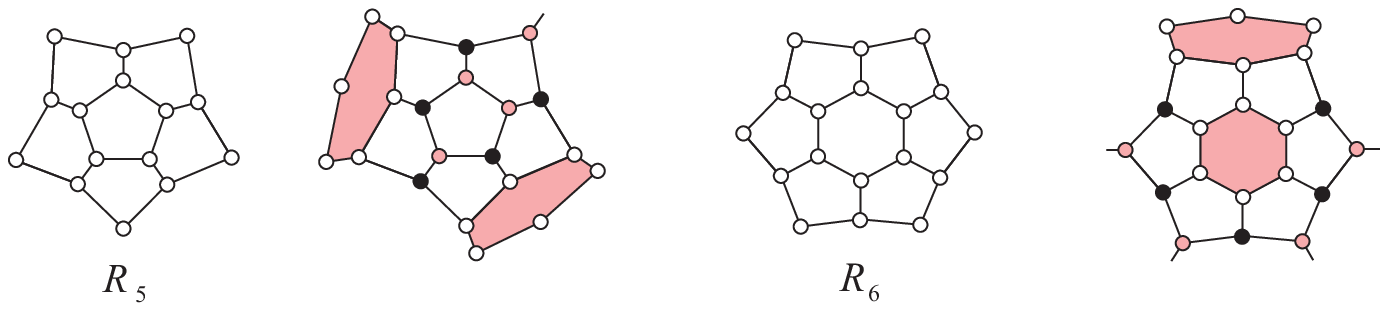}}\\
{Figure \ref{fig3-4}: $R_5$ and $R_6$ and the illustration for the
proof of Theorem  \ref{thm3-3}.}
\end{center}
\end{figure}

On the other hand, we can construct infinitely many fullerene graphs
which are not 2-resonant.
Let $R_5$ and $R_6$ be the graphs obtained by deleting the outer
pentagon from $F_{20}$ and by deleting the outer hexagon from
$F_{24}$, respectively (see Figure \ref{fig3-4}).

\begin{thm}\label{thm3-3}
Let $F$ be a fullerene graph different from $F_{20}$ and $F_{24}$.
If $F$ contains $R_5$ or $R_6$ as subgraphs, then $F$ is not
2-resonant.
\end{thm}

\pf First suppose $R_5\subset F$. Since $F$ is different from
$F_{20}$, there are at least two disjoint hexagons of $F$ adjoining
 $R_5$. Let $\mathcal H$ be the set of
these two hexagons (shadowed hexagons in Figure \ref{fig3-4}). Then
there is a set $S$ of four vertices of Figure \ref{fig3-4}  such
that  $F-\mathcal H-S$ contains five isolated vertices (black
vertices of $R_5$ in Figure \ref{fig3-4}). So $\mathcal H$ is not a
resonant pattern.

Now suppose $R_6\subset F$. Since $F$ is different from $F_{24}$, at
least one hexagon of $F$ adjoins $R_6$. Let $\mathcal H$ be the set
consisting of this hexagon together with the center hexagon of
$R_6$. Similarly, it is easy to see that $\mathcal H$ is not a
resonant pattern (see Figure \ref{fig3-4}). \qed

\vs

Using $R_5$ and $R_6$ as caps, we can construct infinitely many
non-2-resonant nanotubes,  which are, of course,  1-resonant
fullerene graphs. It is interesting to characterize 2-resonant
fullerene graphs. Since each leapfrog fullerene graph is 2-resonant
and has no adjacent pentagons, we now propose  an open problem as
follows.

\begin{prob}
Is every fullerene graph without adjacent pentagons 2-resonant?
\end{prob}


\section{Substructures of $3$-resonant fullerene graphs}

\begin{figure}[!hbtp]\refstepcounter{figure}\label{fig4-1}
\begin{center}
\includegraphics{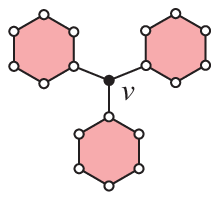}\\
{Figure \ref{fig4-1}: A forbidden subgraph $G^*$ of  3-resonant
fullerene graphs.}
\end{center}
\end{figure}

We first present a  forbidden subgraph  $G^*$ as shown in Figure
\ref{fig4-1} of 3-resonant  fullerene graphs: The three hexagons of
$G^*$ are not  mutually resonant since deleting the three hexagons
isolates the vertex $v$. Let $f$ be a face of a fullerene graph $F$.
A vertex $v$ outside $f$ is {\em adjacent} to $f$ if $v$ has a
neighbor (a vertex adjacent to $v$) in the boundary of $f$. Hence
the forbidden subgraph can be described a vertex being adjacent to
each of  three disjoint hexagons.

\begin{thm}\label{thm4-1}
Let $F$ be a 3-resonant  fullerene graph. Then $|V(F)|\le 60$.
\end{thm}

\pf Since $F$ is 3-resonant, then $F$ contains no $G^*$.  So any
$v\in V(F)$ is adjacent to at least one pentagon of $F$. On the
other hand, for any pentagon $f$ of $F$, there are at most 5
vertices in $V(F-V(f))$ adjacent to it. Hence $|V(F)|\le 12\times
5=60$ since $F$ has exactly 12 pentagons. So the theorem holds.
\qed\\

We now discuss maximal pentagonal fragments and pentagonal rings as
substructures of fullerene graphs in next two subsections, which
will play important roles in construction of 3-resonant fullerene
graphs.

\subsection{Pentagonal fragments}
A {\em fragment $B$} of a fullerene graph $F$ is a
subgraph of $F$ consisting of a cycle together with its interior.
A fragment $B$ is said  to be {\em pentagonal } if its every inner
face is a pentagon. A pentagonal fragment $B$ of a fullerene graph
$F$ is  {\em maximal} if all faces adjoining $B$ are hexagons. For a
pentagonal fragment $B$, use $\gamma (B)$ denote the minimum number
of pentagons adjoining a pentagon in $B$. For example, $\gamma
(R_5)=3$.

 The following two lemmas due to Ye and Zhang are useful.

\begin{lemma}\cite{YZ}\label{lem4-3}
Let $B$ be a fragment of  a fullerene graph $F$ and $W$  the set of
$2$-degree vertices on the boundary $\partial B$. If $0<|W|\le 4$,
then
$T=F-(V(B)\setminus W)$ is a forest and\\
{\upshape{(1)}} $T$ is $K_2$ if $|W|=2$;\\
{\upshape{(2)}} $T$ is $K_{1,3}$ if $|W|=3$;\\
{\upshape{(3)}} $T$ is the union of two $K_2$'s, or a $3$-length
path, or $T_0$ as shown in Figure 6 if $|W|=4$.
\end{lemma}

\begin{figure}[!hbtp]\refstepcounter{figure}\label{fig4-2}
\begin{center}
\includegraphics{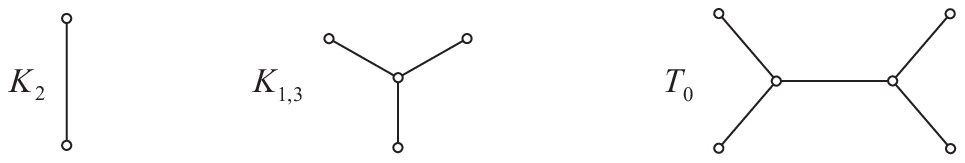}\\
{Figure \ref{fig4-2}: Trees $K_2$, $K_{1,3}$ and $T_0$.}
\end{center}
\end{figure}

\begin{lemma}\cite{YZ}\label{lem4-6}
Let $B$ be a pentagonal fragment of a fullerene graph $F$. Then\\
{\upshape(1)} $R_5\subseteq B$ if $\gamma(B)\ge 3$;\\
{\upshape(2)} $B$ has a pentagon adjoining exactly two adjacent
pentagons of $B$ if $\gamma(B)=2$.
\end{lemma}

\begin{figure}[!hbtp]\refstepcounter{figure}\label{fig4-7}
\begin{center}
\includegraphics{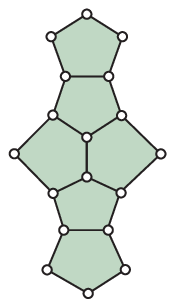}\\
{Figure \ref{fig4-7}: The turtle.}
\end{center}
\end{figure}

A {\em turtle} is a pentagonal fragment consisting of six pentagons
as illustrated in Figure \ref{fig4-7}.  $\gamma(B)=1$ if $B$ is a
turtle. The following theorem characterizes the maximal pentagonal
fragments of 3-resonant  fullerene graphs.

\begin{thm}\label{thm4-7}
Let $F$ be a 3-resonant  fullerene graph different from $F_{20}$ and
$B$ a maximal pentagonal fragment of $F$. Then $B$ is either a
pentagon or a turtle.
\end{thm}

\pf For  a set $\mathcal H$ of at most three disjoint hexagons of
$F$, we have that $\mathcal H$ is a sextet pattern, that is,
$F-\mathcal H$ has a perfect matching, since $F$ is 3-resonant. This
fact will be used repeatedly. Let $B$ be a maximal pentagonal
fragment of $F$. By Theorem \ref{thm3-3}, $B$ contains no $R_5$.
Lemma \ref{lem4-6} implies $\gamma(B)\leq 2$. If $\gamma(B)=0$, then
$B$ is a pentagon. So suppose that $\gamma(B)>0$.

\begin{figure}[!hbtp]\refstepcounter{figure}\label{fig4-8}
\begin{center}
\includegraphics{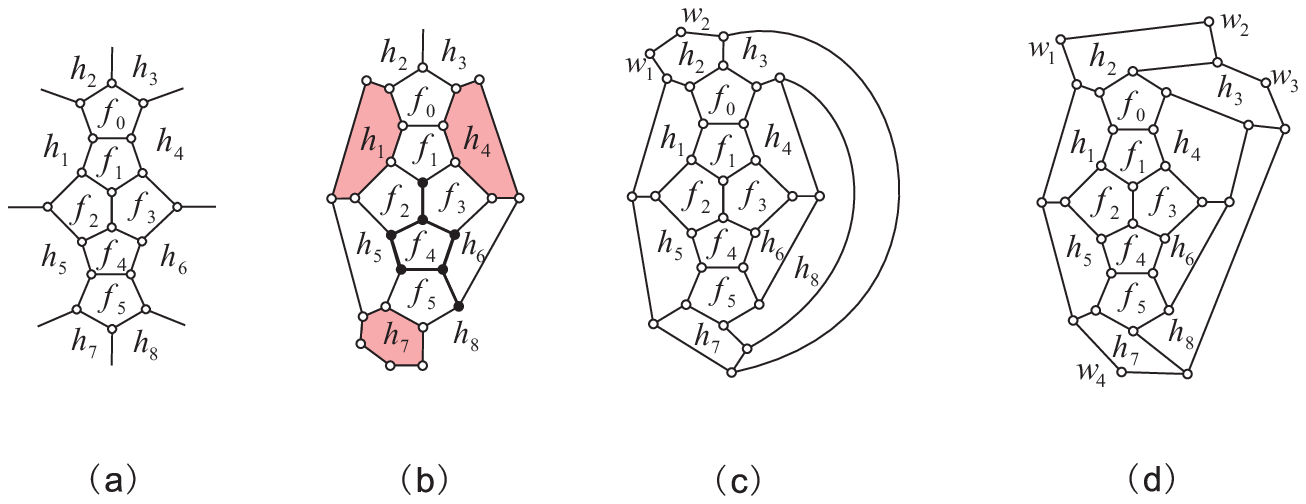}\\
{Figure \ref{fig4-8}: The illustration for the proof of  Case 1 of
Theorem \ref{thm4-7}.}
\end{center}
\end{figure}

{\em Case 1.} $\gamma(B)=1$. Then $B$ has a pentagon $f_0$ with a
unique adjacent pentagon $f_1$. The other four faces adjacent to
$f_0$ are all hexagons since $B$ is maximal, and denoted by
$h_1,h_2,h_3$ and $h_4$ such that $h_i$ is adjacent to $h_{i+1}(1
\leq i\leq 3)$ and both $h_1$ and $h_4$ are also adjacent to $f_1$.
Further, let $f_2$ and $f_3$ be the other faces adjacent to $f_1$ as
illustrated in Figure \ref{fig4-8}(a).

If one of $f_2$ and $f_3$ is a hexagon, say $f_2$, then
$F-\{h_2,h_4,f_2\}$ has an isolated vertex; that is impossible.
Hence both $f_2$ and $f_3$ must be pentagons and thus belong to $B$
since $B$ is maximal. Let $f_4(\ne f_1)$ be the face adjacent to
both $f_2$ and $f_3$. Then  $f_4$ is a pentagon; otherwise,
 $F-\{h_1,h_4,f_4\}$ would have  an isolated vertex.  Let
$f_5$ be the face adjacent to $f_4$ but not adjacent to $f_2$ and
$f_3$. Then $f_5$ is also a pentagon; otherwise,  one component of
$F-\{h_1,h_4,f_5\}$ would be $K_{1,3}$, which has no perfect
matchings. Thus $G:=\cup_{i=0}^5f_i\subseteq B$ is a turtle. It
suffices  to show that $B=G$; that is, all faces adjoining  $G$ are
hexagons.
Besides the four faces $h_1,\dots,h_4$, let $h_5,h_6,h_7$ and $h_8$
be the remaining four  faces adjoining $G$ as illustrated in Figure
\ref{fig4-8}(a). It can be seen that $h_1, \dots, h_8$ are different
from each other. Since $h_1, h_2, h_3$ and $h_4$ are
 hexagons, it remains to show that $h_5, h_6, h_7$ and
$h_8$ are hexagons. Let $G':=G\cup (\cup_{i=1}^8 h_i)$.

We claim that  both $h_5$ and $h_6$ are hexagons. Since
$R_5\nsubseteq B$, one of $h_5$ and $h_6$ must be a hexagon, say
$h_5$, by the symmetry of $G$. Suppose to the contrary that $h_6$ is
a pentagon. Then   $h_7$ is a pentagon; otherwise, $h_1, h_4$ and $
h_7$ are disjoint hexagons, and $F-\{h_1, h_4, h_7\}$ would have an
odd component with seven vertices (see Figure \ref{fig4-8}(b)). If
$h_8$ is a pentagon, then $G'$ is a fragment with only two 2-degree
vertices $w_1$ and $w_2$ on $h_2$ (see Figure \ref{fig4-8}(c)). This
contradicts that $F$ is 3-edge connected since the two edges coming
 out $G_1$ from $w_1$ and $w_2$ form an edge-cut of $F$.  So $h_8$
must be a hexagon and the fragment $G'$ contains four 2-degree
vertices $w_1,w_2,w_3$ and $w_4$ (see Figure \ref{fig4-8}(d)). By
Lemma 4.2(3), there are at most four faces of $F$ outside $G'$.
These faces must be all pentagons since the fragment $G'$ contains
exactly eight pentagons. So $w_1$ and $w_4$ must be adjacent in $F$,
and  $w_2$ and $w_3$ are also adjacent in $F$, resulting in a face
of $F$ with size three. This contradiction establishes the claim.

\begin{figure}[!hbtp]\refstepcounter{figure}\label{fig4-9}
\begin{center}
\includegraphics{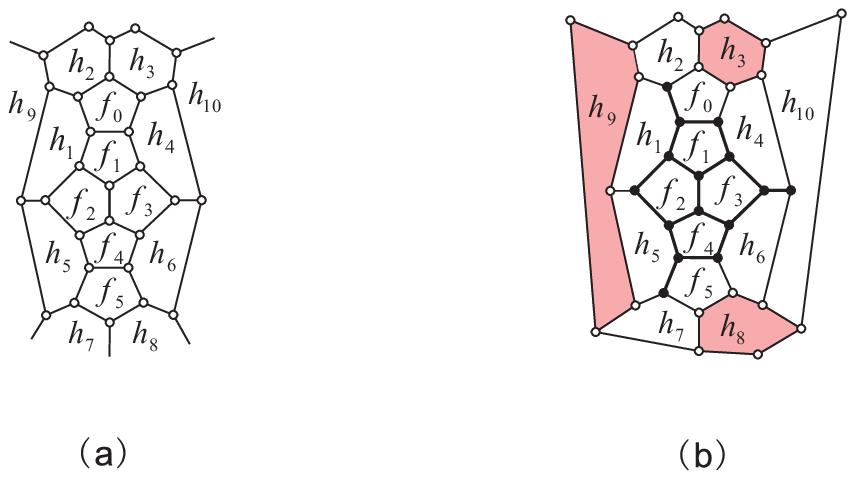}\\
{Figure \ref{fig4-9}: The illustration for the proof of Case 1 of
Theorem \ref{thm4-7}.}
\end{center}
\end{figure}
Further, we claim that both $h_7$ and $h_8$ are hexagons. Without
loss of generality, suppose to the contrary that $h_7$ is a
pentagon. The faces $h_9$ and $h_{10}$  faces of $F$ adjoining $G'$
as shown in Figure \ref{fig4-9}(a) are distinct and disjoint. Then
 $G'':=G'\cup h_9\cup h_{10}$ is a fragment.  If both $h_8$ and $h_9$ are
hexagons, then $ h_3,h_8$ and $h_9$ are disjoint by Lemma
\ref{lem2-3}, and $F-\{h_3,h_8,h_9\}$ would have an odd component
with 15 vertices
 (see Figure \ref{fig4-9}(b)).
Hence at least one of $h_8$ and $h_9$  is a pentagon, and $G''$ is a
fragment with at most four and at least two 2-degree vertices. By
Lemma 4.2, it can be analyzed  analogously that that $G''$ can not
be a subgraph of $F$.   Hence both $h_7$ and $h_8$ are hexagons. So
all faces of $F$ adjoining  $G$ are hexagons and $B=G$.

\begin{figure}[!hbtp]\refstepcounter{figure}\label{fig4-10}
\begin{center}
\includegraphics{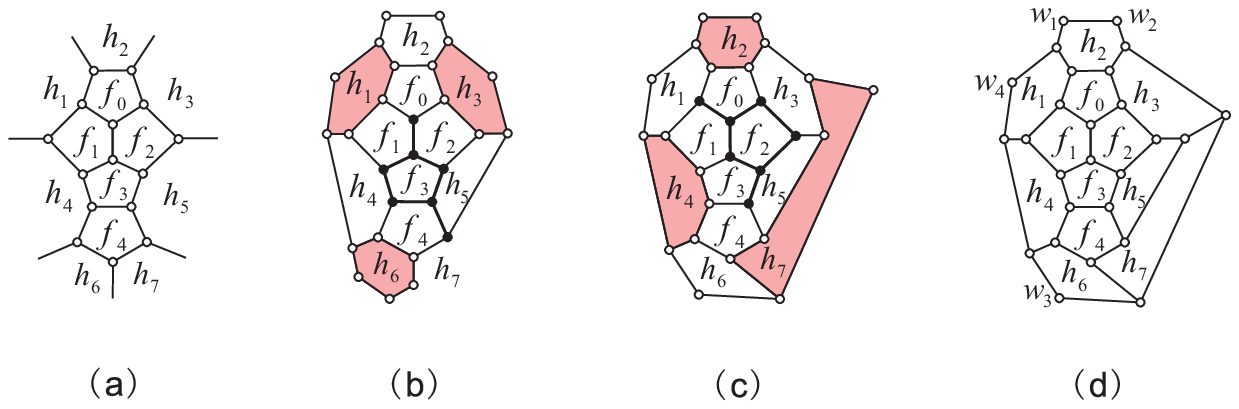}\\
{Figure \ref{fig4-10}: The illustration for the proof of  Case 2 of
Theorem \ref{thm4-7}.}
\end{center}
\end{figure}

{\em Case 2.} $\gamma(B)=2$. Lemma \ref{lem4-6} implies that $B$
contains a pentagon $f_0$ which has exactly two adjacent pentagons
$f_1$ and $f_2$ in $B$. Let $h_1, h_2$ and $h_3$ be the other faces
(hexagons) adjacent to $f_0$ as shown in Figure \ref{fig4-10}(a).

Let $f_3 (\ne f_0)$ be the face of $F$ adjacent to both $f_1$ and
$f_2$. Similarly,  $f_3$ is a pentagon; otherwise,
 disjoint hexagons $h_1,h_3$ and $f_3$ are not mutually resonant.
   Let $h_4,f_4,h_5$ be the other adjacent faces of
$f_3$ as shown in Figure \ref{fig4-10}(a). If $f_4$ is a hexagon,
then  one component of $F-\{h_1,h_3,f_4\}$ is $K_{1,3}$. So $f_4$ is
also a pentagon in $B$. Since $R_5\nsubseteq B$, at least one of
$h_4$ and $h_5$ is a hexagon, say $h_4$. If $h_5$ is also a hexagon,
then  one component of $F-\{h_2, h_4, h_5\}$ is $K_{1,3}$. So $h_5$
is a pentagon.

Let $h_6$ and $h_7$ be the other two adjacent faces of $f_4$ as
shown in Figrue \ref{fig4-10}. If $h_6$ is a hexagon, then $\{h_1,
h_3, h_6\}$ is not a resonant pattern since $F-\{h_1, h_3, h_6\}$
has an odd component with seven vertices (see Figure
\ref{fig4-10}(b)). Hence $h_6$ is a pentagon. Similarly, $h_7$ must
 be a pentagon; if not,  $\{h_2, h_4, h_7\}$ is not a resonant
 pattern (see Figure \ref{fig4-10}(c)). Now, we have a
fragment $G:=(\cup_{i=0}^4f_i)\cup (\cup_{j=1}^7h_j)$ with four
2-degree vertices $w_1,w_2,w_3$ and $w_4$ (see Figure
\ref{fig4-10}(d)). By Lemma \ref{lem4-3}(3), it can be similarly
checked that $G\nsubseteq F$; that is, $\gamma(B)=2$ is impossible.
\qed

\subsection{Pentagonal rings}

\begin{figure}[!hbtp]\refstepcounter{figure}\label{fig4-11}
\begin{center}
\includegraphics{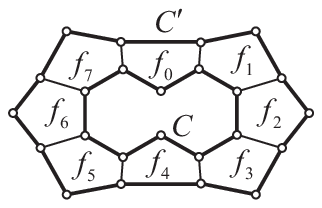}\\
{Figure \ref{fig4-11}: A pentagonal ring $R$ of length  eight with
$s(R)=2$ and $s'(R)=6$.}
\end{center}
\end{figure}

For an integer $l\geq 3$, let $\{f_i| i\in \mathbb{Z}_l\}$ be a
cyclic sequence of $l$ faces (polygons) of a fullerene graph $F$
such that two consecutive faces $f_i$ and $ f_{i+1}$ ($i\in
\mathbb{Z}_l$) intersect only  at an edge, denoted by $e_i$, and two
non-consecutive faces $f_i$ and $f_j$ are disjoint. The subgraph
$R:=\cup_{i\in \mathbb Z_l} f_i$ is called a {\em polygonal ring} of
$F$ if $\{e_i|i\in \mathbb{Z}_l\}$ is a matching of $F$, and $l$ is
called the {\em length} of the polygonal ring $R$, denoted by
$l(R)$. A polygonal ring $R$ is called a {\em pentagonal ring} if
every $f_i$ of $R$  is a pentagon ($i\in \mathbb Z_{l(R)}$) (see
Figure \ref{fig4-11}). The $R_5$ and $R_6$ in Figure \ref{fig3-4}
are two pentagonal rings with length five and six, respectively.

Let $R$ be a pentagonal ring of $F$ consisting of pentagons
$f_1,\dots,f_{l(R)}$. As a subgraph of $F$, $R$ has two faces
different from the $f_i$ ($i=1,\dots,l(R)$). Without loss of
generality, we suppose that $C$ and $C'$ are the boundaries of the
central interior face and exterior face, respectively, and $C$ and
$C'$ have $s(R)$ and $s'(R)$ 2-degree vertices, respectively, with
$s(R)\leq s'(R)$. We call $C$ and $C'$ the {\em inner cycle} and the
{\em outer cycle} of $R$, respectively. Then $s'(R)+s(R)=l(R)$,
$s(R)\leq \lfloor \frac{l(R)}{2}\rfloor$, and $s(R)\ne 1$ and
$s'(R)\ne 1$.

Let $G$ be the subgraph of $F$ induced by the vertices on $C$ and
its interior, and  $r(R), n_6(R)$ and $n_5(R)$  the numbers of
vertices, hexagons and pentagons within $C$, respectively.

We claim that $r(R)$ and $s(R)$ have the same parity. We have
$$|V(G)|=r(R)+l(R)+s(R),$$ and
$$|E(G)|=\frac{2l(R)+3s(R)+3r(R)}{2}.$$
By Euler's formula $|V(G)|-|E(G)|+|F(G)|=1$, where $
|F(G)|(=n_5(R)+n_6(R))$ is the number of the interior faces of $G$,
we have
$$n_5(R)+n_6(R)=\frac{1}{2}(s(R)+r(R)+2). \eqno (1) $$
Further, by $|E(G)|=\frac{1}{2}(5n_5(R)+6n_6(R)+s(R)+l(R))$, we have
$$5n_5(R)+6n_6(R)=2s(R)+3r(R)+l(R)). \eqno (2) $$
Combining Eqs. (1) and (2), we have  that
$$n_5(R)=6+s(R)-l(R),  \eqno (3)$$  and
$$n_6(R)=l(R)+\frac{1}{2}(r(R)-s(R))-5. \eqno (4)$$
Equation (4) implies that $r(R)\equiv s(R)$ (mod 2).

For a fullerene graph $F$, let
$$\psi_l(F):=\min\{s(R)|\mbox{ $R$
is a pentagonal ring of }F{\mbox{ with length } l}\}. \eqno (5) $$
For example, $\psi_5(F_{20})=0$ and $\psi_6(F_{24})=0$. Further, let

$$\tau (F):=\min\{l(R)| R {\mbox{ is a pentagonal ring of }}
F\}. \eqno (6) $$ For example, $\tau(F_{20})=5$ and
$\tau(F_{24})=6$.

\begin{lemma}\label{tau}
For any fullerene graph $F$ with a pentagonal ring,  $5\leq
\tau(F)\leq 12$.
\end{lemma}

 \pf Because $F$ has exactly 12 pentagons, $\tau(F)\leq 12$.
 Further, if $F$ contains a pentagonal ring $R$ with $l(R)\leq 4$,
 then $s(R)=s'(R)=2$ since $F$ has no squares as faces. Hence
 $l(R)=4$, and by Lemma \ref{lem4-3} (1) $F$ has  two edges connecting
the two 2-degree vertices of $R$ lying on the inner cycle and lying
on the outer
 cycle respectively,  which would result in  one face of size at most four in $F$, a contradiction.
Hence $\tau(F)\geq 5$.  \qed\\

The following lemma is due to Kutnar and Maru\v{s}i\v{c}.

\begin{lemma}\cite{KD}\label{lem5-6}
Let $F$ be a  fullerene graph containing a polygonal ring $R$ of
length five, and let $C$ and $C'$ be the inner cycle and the outer
cycle of $R$,
respectively. Then either\\
\noindent{\upshape (1)} $C$ or $C'$ is the boundary of a face, or\\
\noindent{\upshape (2)} both $C$ and $C'$ are of length 10, and the
five faces of $R$ are all hexagonal.
\end{lemma}

By Lemma \ref{lem5-6} we immediately have

\begin{cor}\label{cor5-3}
If a fullerene graph $F$ contains a pentagonal ring $R$ of length
five, then $R$ is just $R_5$.
\end{cor}

\begin{lemma}\label{lem5-3}
There is no fullerene graph $F$ with $\tau (F)=7$.
\end{lemma}

\begin{figure}[!hbtp]\refstepcounter{figure}\label{fig5-6}
\begin{center}
\includegraphics{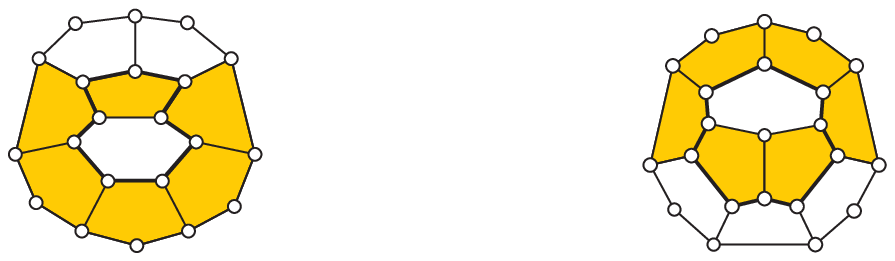}\\
{Figure \ref{fig5-6}: The illustration for the proof of Lemma
\ref{lem5-3}.}
\end{center}
\end{figure}

\pf Suppose to the contrary that $F$ is a fullerene graph with
$\tau(F)=7$. Let $R$ be a pentagonal ring of $F$ with length
$l(R)=7$. Then $s(R)\le \lfloor \frac {l(R)} 2\rfloor=3$. So
$s(R)=2$ or $3$. By Lemma \ref{lem4-3}, whenever $s(R)=2$ or 3, $F$
would contain a $R_6$ (see  Figure \ref{fig5-6}), contradicting that
 $\tau(F)=7$.\qed

\begin{lemma}\label{lem5-8}
A fullerene graph $F$ with $\tau(F)=11$ is not $3$-resonant.
\end{lemma}

\pf Let $R$ be a pentagonal ring of $F$ with length
$l(R)=\tau(F)=11$ and $s(R)=\psi_{11}(F)$. Let $C$ be the inner
cycle of $R$. By Eq. (3), $n_5(R)=s(R)-5$, and $\psi_{11}(F)=s(R)\ge
5$. On the other hand, $\psi_{11}(F)=s(R)\leq \lfloor
\frac{l(R)}{2}\rfloor =5$. So $\psi_{11}(F)=s(R)=5$ and $n_5(R)=0$;
that is, there are no pentagons within $C$.

Let $v_1, v_2, v_3, v_4$ and $v_5$ be the five 2-degree vertices
clockwise on $C$. If two of these five vertices  are adjacent in
$F$, then by Lemma \ref{lem4-3} it follows that the two vertices are
consecutive, say $v_1$, $v_2$, and the other three vertices
$v_3,v_4$ and $v_5$ have a common neighbor within $C$, denoted by
$w$. Let $h$ be the face of $F$ containing $v_1,v_2,v_3,w$ and
$v_5$. Note that any two of $v_1,v_2,\dots, v_5$ are not adjacent on
$C$. So $|h|\ge 7$, a contradiction. If any two of $v_1,v_2,\dots,
v_5$ have no common neighbor within $C$, then the five faces
adjoining $R$ along $C$ form a polygonal ring $R'$ with $C$ as the
outer cycle. Since $|C|=16$, the inner cycle of $R'$ bounds a face
$f'$ of $F$ by Lemma \ref{lem5-6}. Note $s(R)\equiv r(R)$ (mod 2).
So $f'$ is a pentagon, contradicting $n_5(R)=0$.

\begin{figure}[!hbtp]\refstepcounter{figure}\label{fig5-14}
\begin{center}
\includegraphics{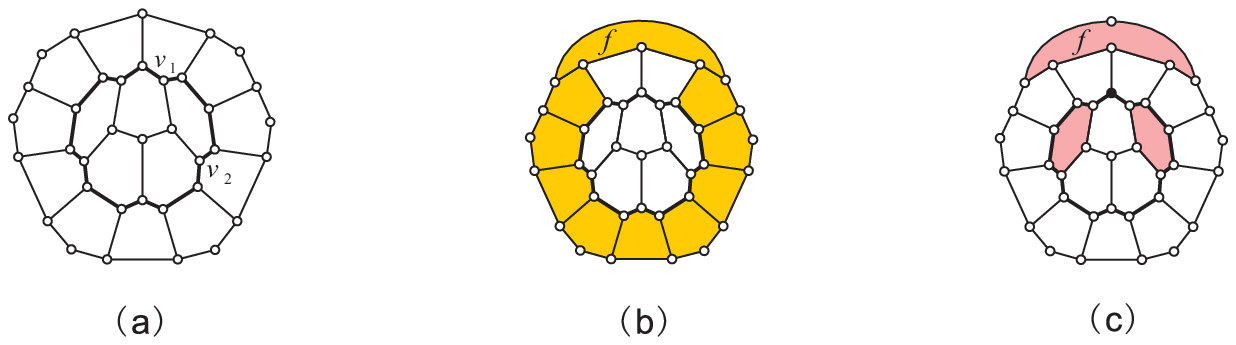}\\
{Figure \ref{fig5-14}: The illustration of Lemma \ref{lem5-8}.}
\end{center}
\end{figure}

So there exist two vertices of $v_1,\dots,v_5$ with a common
neighbor within $C$.  They must be consecutive by Lemma
\ref{lem4-3}, so say $v_1$ and $v_2$. By Lemma \ref{lem4-3} and
$n_5(R)=0$, the subgraph of $F$ induced by $R$ together with all
vertices within $C$ is isomorphic to the graph  in Figure
\ref{fig5-14}(a). Let $f$ be the  face adjacent to $R$ along a
4-length path on the boundary of $R$ (see Figure \ref{fig5-14}(b)).
If $f$ is a pentagon, then $F$ contains a pentagonal ring $R'$ with
length $l(R')=10$ (see Figure \ref{fig5-14}(b)). Then $11=\tau
(F)\le l(R')=10$, that is a contradiction. So suppose $f$ is a
hexagon. Then $F$ contains the forbidden subgraph of 3-resonant
fullerene graph in Figure \ref{fig4-1}; see also  Figure
\ref{fig5-14}(c). Hence $F$ is not 3-resonant.  \qed


\section{Construction of $k$-resonant ($k\ge 3$) fullerene graphs}

For a pentagon $f$ of a fullerene graph $F$, if it dose not lie in
any pentagonal ring of $F$, then it must lie in some maximal
pentagonal fragment of $F$. In particular, if $F$ is a 3-resonant
fullerene graph containing no pentagonal rings, then by  Theorem
\ref{thm4-7} the  maximal pentagonal fragment  of $F$ containing any
given pentagon is either a pentagon or a turtle.

\begin{figure}[!hbtp]\refstepcounter{figure}\label{fig5-1}
\begin{center}
\scalebox{0.8}{\includegraphics{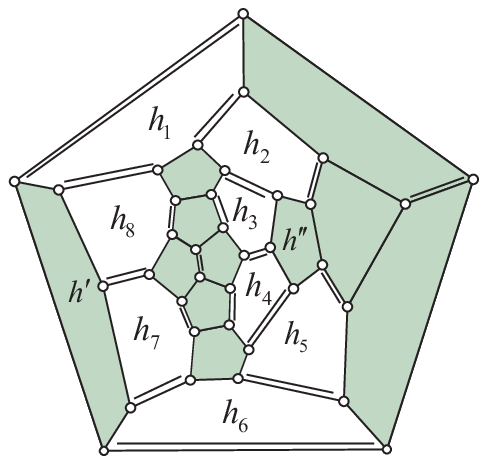}}\\
{Figure \ref{fig5-1}: The fullerene graph $F_{36}^1$ with a perfect
matching $M$.}
\end{center}
\end{figure}

\begin{lemma}\label{lem5-1}
Let $F$ be a  fullerene graph without pentagonal rings. Then $F$ is
$3$-resonant if and only if $F$ is either $\text{C}_{60}$ or
$F_{36}^1$ shown in Figure \ref{fig5-1}. Further, both
$\text{C}_{60}$ and  $F_{36}^1$ are $k$-resonant for any integer
$k\geq 3$.
\end{lemma}

\pf Let $F$ be a 3-resonant  fullerene graph without pentagonal
rings as subgraphs. Then $F$ is different from $F_{20}$ since
$F_{20}$ contains a pentagonal ring $R_5$. So by Theorem
\ref{thm4-7}, every maximal pentagonal fragment of $F$ is  either a
pentagon or a turtle. If $F$ contains no turtles as maximal
pentagonal fragments, then every pentagon of $F$ is adjacent only to
hexagons. Hence $F$ satisfies IPR (isolated pentagon rule). By
Theorem \ref{thm4-1}, $F$ is $\text{C}_{60}$ since it is the unique
fullerene graph with no more than 60 vertices and without adjacent
pentagons.

\begin{figure}[!hbtp]\refstepcounter{figure}\label{fig5-2}
\begin{center}
\includegraphics{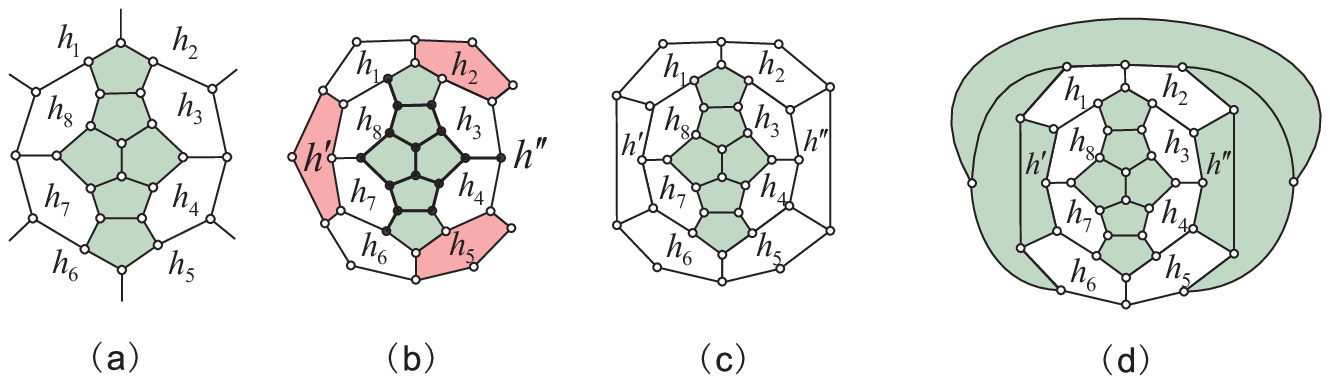}\\
{Figure \ref{fig5-2}: The illustration for the proof of Lemma
\ref{lem5-1}.}
\end{center}
\end{figure}

Now suppose that $F$ contains a turtle $B$ as a maximal pentagonal
fragment. Denote clockwise the hexagons adjoining  $B$ by $h_1, h_2,
..., h_8$ as shown in Figure \ref{fig5-2}(a). Let $G_0:=B\cup
h_3\cup h_4\cup h_7\cup h_8$. Then $h_1,h_2, h_5$ and $h_6$ are four
hexagons adjoining  $G_0$. The other two faces adjoining  $G_0$ are
denoted by $h'$ and $h''$ such that $h'$ is adjacent to both $h_7$
and $h_8$. By Lemma $\ref{lem2-3}$, $h'$ is disjoint from $h_2$ and
$h_5$. If $h'$ is a hexagon, then $\mathcal H=\{h', h_2, h_5\}$ is
not a resonant pattern since $F-\mathcal H$ has a component with 15
vertices (see Figure \ref{fig5-2}(b)), contradicting  that $F$ is
$3$-resonant. So $h'$ must be a pentagon. By the symmetry of $G_0$,
$h''$ is also pentagonal. Hence the fragment $G_1$, consisting of
$G_0$ together with its all adjacent faces, has exactly four
2-degree vertices on its boundary (see Figure \ref{fig5-2}(c)). By
Lemma \ref{lem4-3}(3), $F$ is isomorphic to the graph (d) in Figure
\ref{fig5-2}, that is  $F_{36}^1$ in Figure \ref{fig5-1}.

Conversely,  each of fullerene graphs $\text{C}_{60}$ and $F_{36}^1$
has a perfect matching, illustrated by double edges in Figures
\ref{fig3-1} and \ref{fig5-1} respectively, so that all hexagons are
alternating. Hence $\text{C}_{60}$ and $F_{36}^1$ are $k$-resonant
for any integer $k\geq 1$ since any disjoint hexagons are mutually
resonant. \qed

From now on we discuss $3$-resonant fullerene graphs with a
pentagonal ring. By Lemmas \ref{tau}, \ref{lem5-3} and \ref{lem5-8},
we have that $\tau(F)=5,6,8,9,10$ or 12. Such cases will be
discussed in next five lemmas.

\begin{lemma}\label{lem5-2}
Let $F$ be a  fullerene graph with $\tau(F)=5$ or 6. Then $F$ is
$3$-resonant if and only if it is either $F_{20}$ or $F_{24}$
(Figure \ref{fig3-2}). Further, $F_{20}$ and $F_{24}$ are
$k$-resonant for any integer $k\geq 3$.
\end{lemma}

\pf Since both $F_{20}$ and $F_{24}$ are 2-resonant and contain no
more than two hexagons, they are also $k$-resonant for any integer
$k\geq 3$.

Now let $F$ be a  $3$-resonant  fullerene graph. If $\tau (F)=5$,
then $F$ contains pentagonal ring $R_5$ by Corollary \ref{cor5-3}.
So $F$ is  $F_{20}$ by Theorem \ref{thm3-3}.

\begin{figure}[!hbtp]\refstepcounter{figure}\label{fig5-5}
\begin{center}
\includegraphics{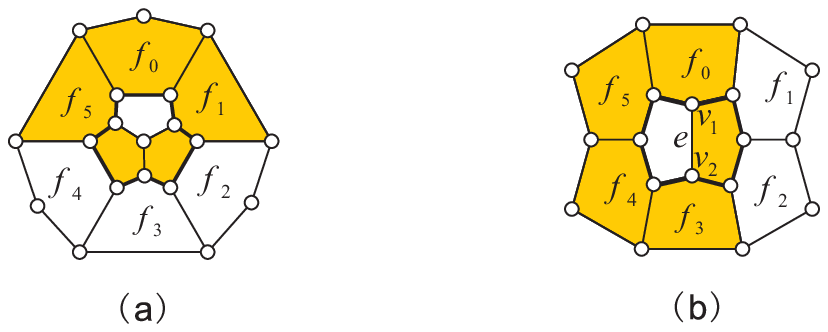}\\
{Figure \ref{fig5-5}: The illustration for the proof of Lemma
\ref{lem5-2}.}
\end{center}
\end{figure}
Now suppose $\tau (F)=6$. Let $R$ be a pentagonal ring with length
$l(R)=\tau(F)=6$ and let $C$ and $C'$ be the inner cycle and the
outer cycle of $R$, respectively. Let $f_0,f_1,\dots, f_5$ be the
six pentagons of $R$ in clockwise order. Then $1\ne s(R)\leq \lfloor
\frac{6}{2}\rfloor=3$. If $s(R)=0$, then $R$ is  $R_6$ and $F$ is
just  $F_{24}$ by Theorem \ref{thm3-3}.

If $s(R)=3$, there are three 2-degree vertices on $C$ and also three
2-degree vertices on $C'$. By Lemma \ref{lem4-3}, the three 2-degree
vertices on $C$ have a common neighbor within $C$. Hence $F$
contains a $R_5$ (see Figure \ref{fig5-5}(a)),  contradicting
$\tau(F)=6$. If $s(R)=2$, there are two 2-degree vertices $v_1,v_2$
on $C$. By Lemma \ref{lem4-3}, $v_1$ and $v_2$ are adjacent in $F$
(see Figure \ref{fig5-5}(b)). Hence  $F$ contains a $R_5$, also
contradicting $\tau (F)=6$. \qed

\begin{figure}[!hbtp]\refstepcounter{figure}\label{fig5-19}
\begin{center}
\includegraphics{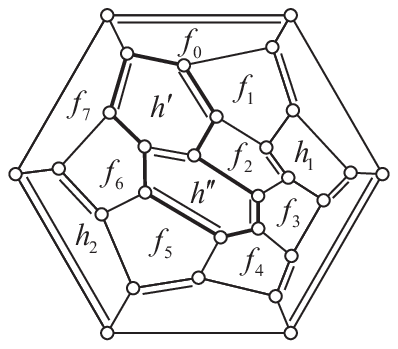}\\
{Figure \ref{fig5-19}: The fullerene graph $F_{28}$ with a perfect
matching.}
\end{center}
\end{figure}

\begin{lemma}\label{lem5-4}
Let $F$ be a fullerene graph with $\tau(F)=8$. Then  $F$ is
$3$-resonant if and only if $F$ is   $F_{28}$ shown in Figure
\ref{fig5-19}. Further, $F_{28}$ is also $k$-resonant for any
integer $k\geq 3$.
\end{lemma}

\pf Similar to the proof of Lemma \ref{lem5-1} we can show readily
that $F_{28}$ is $k$-resonant for any integer $k\geq 3$.

Conversely, let  $F$  be a $3$-resonant fullerene graph with
$\tau(F)=8$. Let $R$ be a pentagonal ring of $F$ with length
$l(R)=\tau(F)=8$ and $s(R)=\psi_8(F)$. Let $C$ and $C'$ be the inner
cycle and the outer cycle of $R$, respectively, and $f_0, f_1,
\dots, f_7$ the eight pentagons of $R$ in clockwise order.
Obviously, $2\leq \psi_8(F)=s(R)\leq \lfloor
\frac{l(R)}{2}\rfloor=4$.

\begin{figure}[!hbtp]\refstepcounter{figure}\label{fig5-9}
\begin{center}
\includegraphics{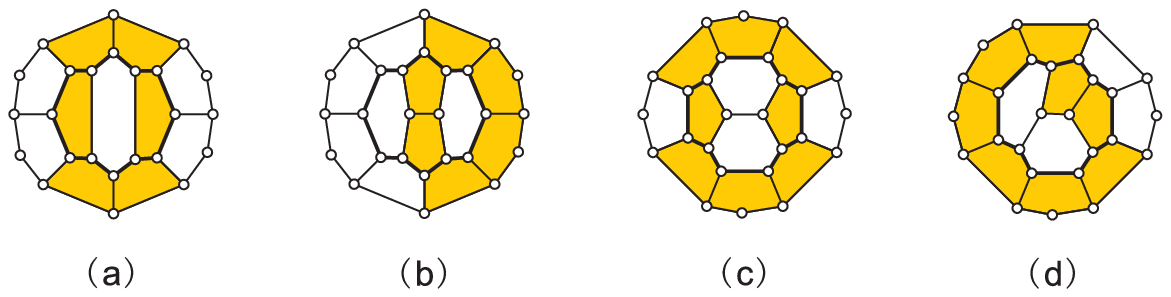}\\
{Figure \ref{fig5-9}: The illustration for Case 1 of the proof of
Lemma \ref{lem5-4}.}
\end{center}
\end{figure}

{\em Case 1.}  $\psi_8(F)=4$. By Lemma \ref{lem4-3}(3), the subgraph
$G$ of $F$ induced by $R$ together with all vertices within $C$ is
isomorphic to one of the four graphs shown in Figure \ref{fig5-9}.
If $G$ is isomorphic to the graph (a) or (b), then $F$ contains a
pentagonal ring with length six, contradicting $\tau(F)=8$. If $G$
is isomorphic to the graph (c) or (d), then $F$ contains a
pentagonal ring $R'$ with length eight and $s(R')=2$, contradicting
$\psi_8(F)=4$ (refer to Eq. (5)).

\begin{figure}[!hbtp]\refstepcounter{figure}\label{fig5-8}
\begin{center}
\includegraphics{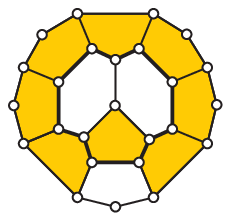}\\
{Figure \ref{fig5-8}: The illustration for Case 2 of the proof of
Lemma \ref{lem5-4}.}
\end{center}
\end{figure}

{\em Case 2.} $\psi_8(F)=3$. By Lemma \ref{lem4-3} (2) and Eqs. (3)
and (4), we have $r(R)=1$, $n_5(R)=1$ and $n_6(R)=2$. Hence $F$
contains a pentagonal ring $R'$ with length eight and $s(R')=2$; see
 Figure \ref{fig5-8}. So
$3=\psi_8(F)\le s(R')$=2 by (5), a contradiction.

\begin{figure}[!hbtp]\refstepcounter{figure}\label{fig5-7}
\begin{center}
\includegraphics{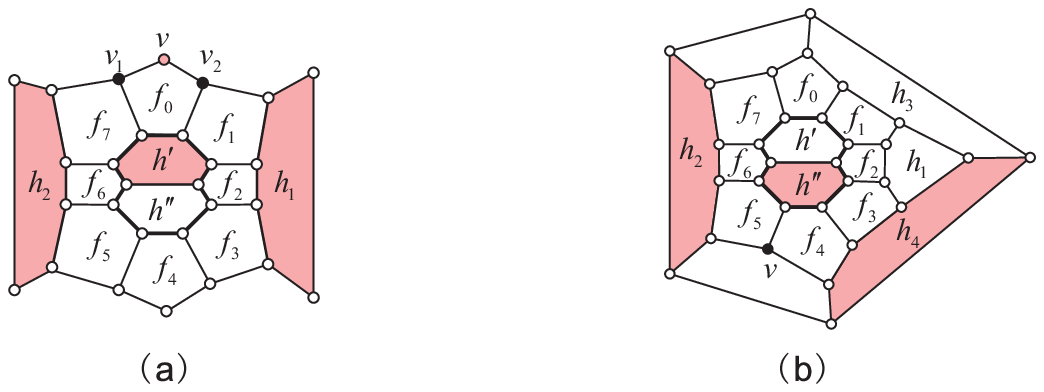}\\
{Figure \ref{fig5-7}: The illustration for Case 3 of the proof of
Lemma \ref{lem5-4}.}
\end{center}
\end{figure}

{\em Case 3.} $\psi_8(F)=2$. Then $R$ contains two 2-degree vertices
$u_1$ and $u_2$ on $C$. By Lemma \ref{lem4-3}(1), $u_1$ and $u_2$
are adjacent in $F$, and  lie on two pentagons $f_{i}$ and $f_{i+4}$
for some $i\in \mathbb{Z}_8$, respectively (say $i=2$, and see
Figure \ref{fig5-7} (a)). So there are exactly two adjacent hexagons
$h'$ and $h''$ within $C$. Let $h_1$ and $h_2$ be the two faces
outside $C'$ such that $h_1$ is adjacent to faces $f_{1}$, $f_{2}$
and $f_{3}$, while $h_2$ is adjacent to faces $f_{5}$, $f_{6}$ and
$f_{7}$ (see Figure \ref{fig5-7} (a)). Then $h_1$ and $h_2$ are
distinct and disjoint. If both $h_1$ and $h_2$ are hexagons, let $v,
v_1$ and $v_2$ be three vertices on $C'$ as shown in Figure
\ref{fig5-7} (a) and let $S:=\{v\}$ and $\mathcal H:=\{h_1, h_2,
h'\}$. Then $F-\mathcal H-S$ has two isolated vertices $v_1$ and
$v_2$. By Lemma \ref{thm2-1}, $\mathcal H$ is not a resonant
pattern, contradicting that $F$ is $3$-resonant. So at least one of
$h_1$ and $h_2$,  say $h_1$, is a pentagon. If $h_2$ is a hexagon,
let $h_3$ and $h_4$ be the other two adjacent faces of $h_1$ as
shown in Figure \ref{fig5-7} (b). By Lemma \ref{lem4-3}(3), it
follows  that both $h_3$ and $h_4$ are hexagons. Hence $F$ is the
fullerene graph $F_{30}$ shown in Figure \ref{fig5-7} (b). Clearly,
$\mathcal H:=\{h_2, h_4, h''\}$ is not  resonant  since $F-\mathcal
H$ has an isolated vertex $v$ (see Figure \ref{fig5-7} (b)), also
contradicting that $F$ is 3-resonant. So both $h_1$ and $h_2$ are
pentagons. By Lemma \ref{lem4-3}, $F$ is
 the fullerene graph $F_{28}$ shown in Figure
\ref{fig5-19}.\qed

\begin{figure}[!hbtp]\refstepcounter{figure}\label{fig5-18}
\begin{center}
\includegraphics{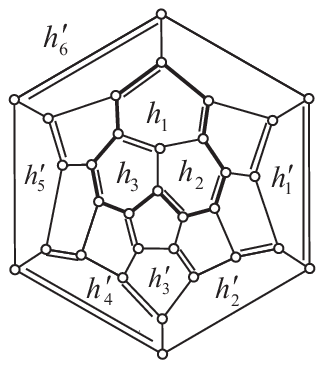}\\
{Figure \ref{fig5-18}: The fullerene graph $F_{32}$ with a perfect
matching.}
\end{center}
\end{figure}

\begin{lemma}\label{lem5-5}
Let $F$ be a fullerene graph with $\tau(F)=9$. Then $F$ is
$3$-resonant if and only if $F$ is  $F_{32}$ as shown in Figure
\ref{fig5-18}. Further, $F_{32}$ is $k$-resonant for any integer
$k\geq 3$.
\end{lemma}

\pf Let  $F$  be a $3$-resonant fullerene graph with $\tau(F)=9$.
Let $R$ be a pentagonal ring of $F$ with length $l(R)=9$ and
$s(R)=\psi_9(F)$. Let $C$ and $C'$ be the inner cycle and the outer
cycle of $R$, respectively. If $s(R)=2$, by by Lemma \ref{lem4-3}
(1)   the two 2-degree vertices on $C$ must be adjacent in $F$. Then
the two faces $f$ and $f'$ within $C$ satisfy $|f|+|f'|=|C|+2=13$.
Hence $F$ has a face of size larger than 6, a contradiction. So
$3\le \psi_9(F)=s(R)\le \lfloor\frac{l(R)} 2\rfloor=4$. Let
$f_0,f_1,\dots,f_8$  be the nine pentagons of $R$ in clockwise
order.

\begin{figure}[!hbtp]\refstepcounter{figure}\label{fig5-11}
\begin{center}
\includegraphics{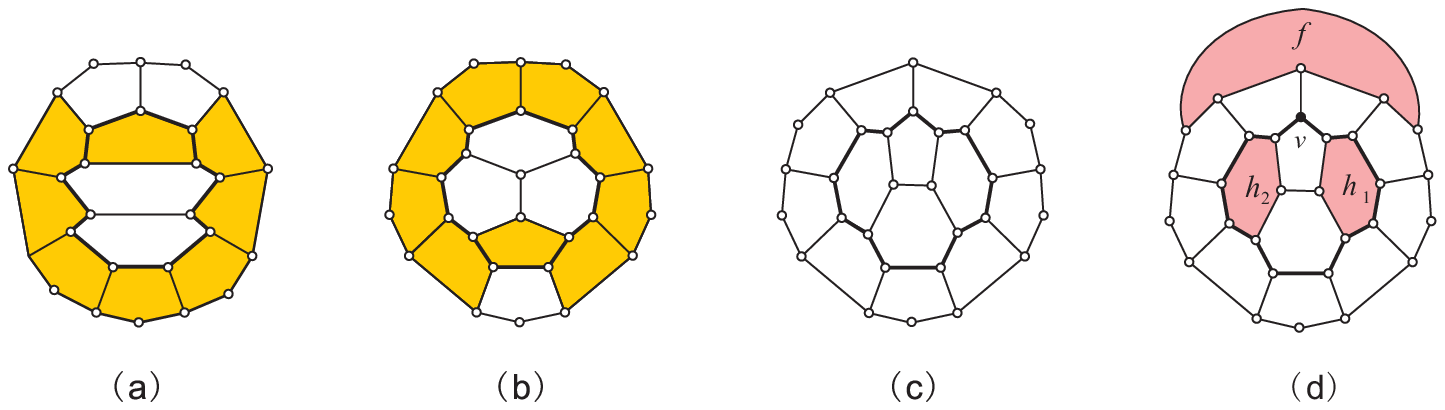}\\
{Figure \ref{fig5-11}: The illustration for Case 1 of the proof
Lemma \ref{lem5-5}.}
\end{center}
\end{figure}

{\em Case 1.} $\psi_9(F)=4$. Let $G$ be the subgraph of $F$ induced
by the vertices of $R$ together with the vertices within $C$. By
Eqs. (3) and (4), $n_5(R)=1$ and $n_6(R)=2+\frac{1}{2}r(R)$. By
Lemma \ref{lem4-3}(3), either $r(R)=0$
 and $n_6(R)=2$, or $r(R)=2$ and $n_6(R)=3$.
By Lemma \ref{lem4-3}, $G$ is isomorphic to the graph (a) in Figure
\ref{fig5-11} if the former holds, and $G$ is isomorphic to either
the graph (b) or (c) in Figure \ref{fig5-11} if the latter holds.

If $G$ is isomorphic to the graph (a), then $F$ contains a pentagon
ring with length eight, contradicting $\tau (F)=9$. If $G$ is
isomorphic to the graph (b), then $F$ contains a pentagonal ring
$R'$ with length 9 and $s(R')=3$, contradicting $\psi_9(F)=4$. So
suppose $G$ is isomorphic to the graph (c). Let $f$ be the face
adjoining $R$ along a 4-length path as shown in Figure \ref{fig5-11}
(d). If $f$ is a pentagon, then $F$ contains a pentagonal ring with
length 8 which consists of seven pentagons of $R$ and $f$,
contradicting $\tau(F)=9$. So $f$ is a hexagon. Then  disjoint
hexagons $f, h_1$ and $h_2$ form a forbidden substructure of
3-resonant fullerene graphs  (see Figures \ref{fig5-11}(d) and
\ref{fig4-1}).

\begin{figure}[!hbtp]\refstepcounter{figure}\label{fig5-10}
\begin{center}
\includegraphics{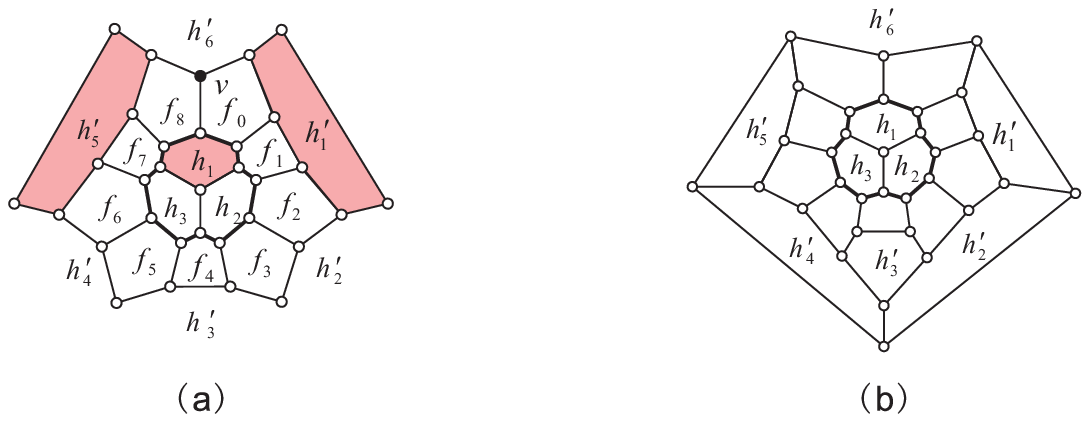}\\
{Figure \ref{fig5-10}: The illustration for Case 2 of the proof of
Lemma \ref{lem5-5}.}
\end{center}
\end{figure}

{\em Case 2.} $\psi_9(F)=3$.  By Lemma \ref{lem4-3}(2) and Eqs (3)
and (4), we have $r(R)=1$, $n_5(R)=0$ and $n_6(R)=3$. Denote the
three hexagons within $C$ by $h_1,h_2$ and $h_3$. The three 2-degree
vertices on $C$ must lie on the pentagons $f_{i}$, $f_{i+3}$ and
$f_{i+6}$ for some $i\in {\mathbb Z}_9$ (say $i=1$, and see Figure
\ref{fig5-10} (a)). Let $h_1', h_2', \dots, h_6'$ be the six faces
adjoining $R$ clockwise along $C'$ such that $h'_1$ is adjacent to
$f_{0}$, $f_{1}$ and $f_{2}$ (see Figure \ref{fig5-10}(a)). If at
least two of $h_1', h_3'$ and $h_5'$ are hexagons, say $h_1'$ and
$h_5'$, then $ h_1',h_5'$ and $h_1$ also form a forbidden
substructure of 3-resonant fullerene graphs (see Figure
\ref{fig5-10}(a)), contradicting  that $F$ is $3$-resonant. If only
one of $h_1'$, $h_3'$ and $h_5'$ is a hexagon, then by Lemma
\ref{lem4-3}(3), $F$ has a face with size seven, a contradiction. So
all of $h_1', h_3'$ and $h_5'$ are pentagons. By Lemma
\ref{lem4-3}(2), $F$ is isomorphic to the graph (b) in Figure
\ref{fig5-10}; that is $F_{32}$ in Figure \ref{fig5-18}.

Conversely, it  needs to show that $F_{32}$ is $k$-resonant for any
integer $k\geq 3$. Since there are no more than two disjoint
hexagons in $F_{32}$, it suffices to show that any two disjoint
hexagons of $F_{32}$ are mutually resonant. By symmetry, we only
consider $\{h_1,h'_4\}$ and $\{h_1,h'_6\}$. Let $M$ be the perfect
matching of $F_{32}$ consisting of the double edges illustrated in
Figure \ref{fig5-18}. Clearly,   $h_1,h'_4$ and $h'_6$ are all
$M$-alternating. Hence both $\{h_1,h'_4\}$ and $\{h_1,h'_6\}$ are
resonant patterns of $F_{32}$.\qed

\begin{figure}[!hbtp]\refstepcounter{figure}\label{fig5-20}
\begin{center}
\includegraphics[totalheight=4cm]{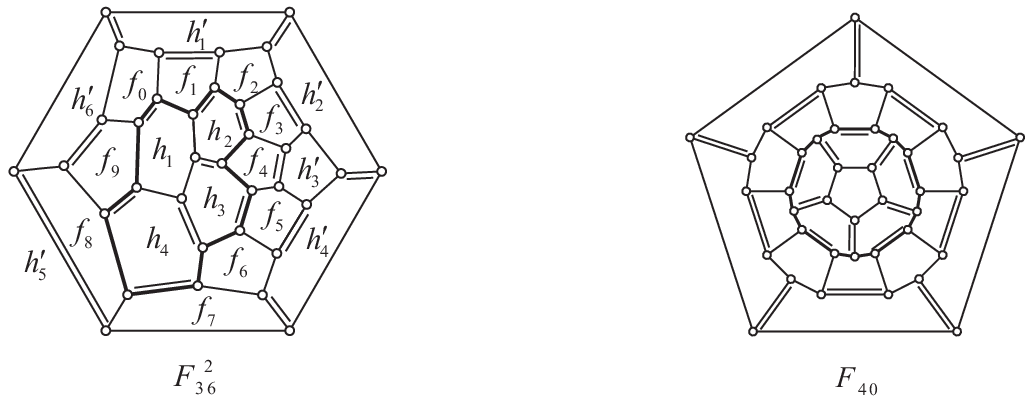}\\
{Figure \ref{fig5-20}: The fullerene graphs $F_{36}^2$ and
$F_{40}$.}
\end{center}
\end{figure}

\begin{lemma}\label{lem5-7}
Let $F$ be a  fullerene graph with $\tau(F)=10$. Then $F$ is
$3$-resonant if and only if $F$ is  either $F_{36}^2$ or $F_{40}$
shown in Figure \ref{fig5-20}. Further, $F_{36}^2$ and $F_{40}$ are
$k$-resonant for any integer $k\geq 3$.
\end{lemma}

\pf Let $F$ be a 3-resonant fullerene graph with $\tau(F)=10$. Let
$R$ be a pentagonal ring of $F$ with length $l(R)=\tau(F)=10$ and
$s(R)=\psi_{10}(F)$. Let $C$ be the inner cycle of $R$ and
$f_0,\dots, f_9$ the pentagons of $R$ in clockwise order. Clearly,
$\psi_{10}(F)\ge 2$. If $\psi_{10}(F)=2$, then the two 2-degree
vertices of $R$ on $C$ are adjacent in $F$ by Lemma \ref{lem4-3}.
Then there are two faces $h_0$ and $h_1$ of $F$ within $C$ such that
 $|h_0|+|h_1|=l(R)+s(R)+2=14$ since the edge within $C$ is counted
twice in $|h_0|+|h_1|$. So at least one of $h_0$ and $h_1$ has size
more than six, a contradiction. If $\psi_{10}(F)=3$, then the three
2-degree vertices of $R$ on $C$ together with vertices within $C$
induce a $K_{1,3}$ by Lemma \ref{lem4-3}(2). Hence there are three
faces $h_0,h_1,h_2$ of $F$ within $C$ and $\sum_{i\in \mathbb
Z_3}|h_i|=l(R)+s(R)+6=19$. So at least one of $h_0$, $h_1$ and $h_2$
has size no less than seven, a contradiction, too. So $4\le
\psi_{10}(F)=s(R)\le \lfloor \frac{l(R)} 2 \rfloor =5$.

{\em Case 1.} $\psi_{10}(F)=4$. By Eqs. (3) and (4), $n_5(R)=0$ and
 $n_6(R)=3+\frac{1}{2}r(R)$. By Lemma
\ref{lem4-3}(3),  $r(R)=0$ or 2.

\begin{figure}[!hbtp]\refstepcounter{figure}\label{fig5-12}
\begin{center}
\includegraphics{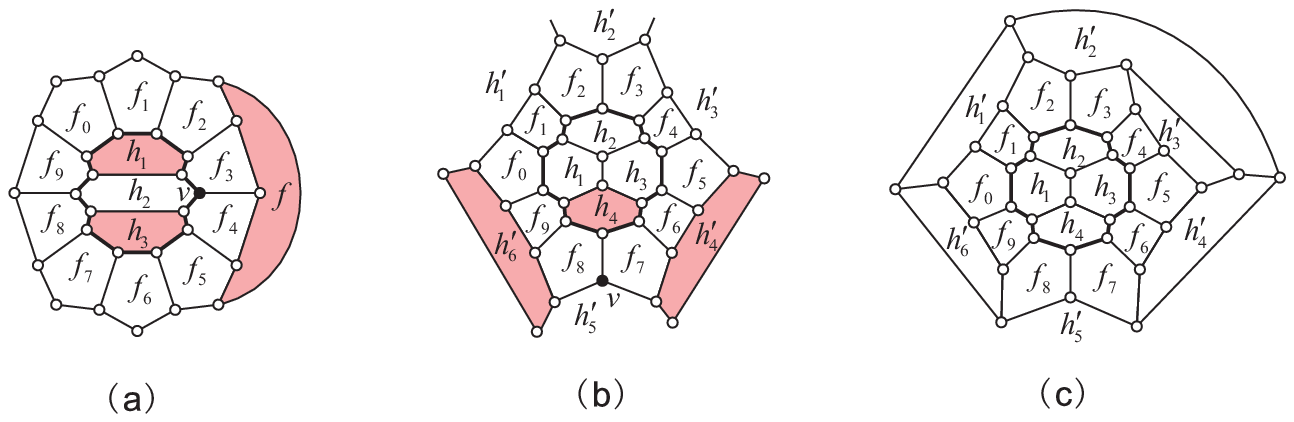}\\
{Figure \ref{fig5-12}: The illustration for Case 1 of the proof of
Lemma \ref{lem5-7}.}
\end{center}
\end{figure}

If $r(R)=0$ and $n_6(R)=3$, the four 2-degree vertices on $C$ belong
to four pentagons $f_{i}, f_{i+1}, f_{i+5}$ and $f_{i+6}$ for some
$i\in \mathbb{Z}_{10}$, say $i=3$ (see Figure \ref{fig5-12}(a)). Let
$h_1, h_2$ and $h_3$ be the three hexagons within $C$ such that
$h_1\cap f_2\ne \emptyset, h_2\cap f_3\ne \emptyset$ and $h_3\cap
f_4\ne \emptyset$. Let $f$ be the common adjacent face  of
$f_2,f_3,f_4$ and $f_5$ (see Figure \ref{fig5-12}(a)). If $f$ is a
pentagon, then $F$ contains a pentagonal ring $(R-\{f_3,f_4\})\cup
f$ with length 9, contradicting $\tau (F)=10$. So suppose $f$ is a
hexagon. Then $\mathcal H:=\{f,h_1,h_3\}$ is not a resonant pattern
since $F-\mathcal H$ has an isolated vertex (the black vertex in
Figure \ref{fig5-12} (a)), contradicting that $F$ is $3$-resonant.

So suppose $r(R)=2$ and $n_6(R)=4$. By Lemma \ref{lem4-3}(3), the
four 2-degree vertices on $C$ together with all vertices within $C$
induce a $T_0$. Hence, the four 2-degree vertices belong to the
pentagons $f_{j}$, $f_{j+3}, f_{j+5}$ and $f_{j+8}$ for some $j\in
\mathbb{Z}_{10}$, say $j=1$ (see Figure \ref{fig5-12}(b)). Let $h_1,
h_2, h_3$ and $h_4$ be the four hexagons within $C$ in clockwise
order and $h_1\cap f_0\ne \emptyset$. Let $h'_1,h'_2,\cdots h'_6$ be
the faces adjoining $R$ in clockwise order along its boundary such
that $h'_1\cap f_1\ne \emptyset$ (see Figure \ref{fig5-12} (b)). By
Lemmas \ref{lem2-3} and \ref{lem5-6}, they are pairwise distinct. If
both $h_4'$ and $h_6'$ are hexagons, then $\mathcal
H:=\{h_4,h'_4,h'_6\}$ is not a resonant pattern since $F-\mathcal H$
has an isolated vertex $v$, contradicting that $F$ is $3$-resonant.
So one of $h'_4$ and $h'_6$ is a pentagon, say $h'_6$. By the
symmetry, one of $h_1'$ and $h'_3$ is a pentagon. If $h'_1$ is a
pentagon, then $F$ would have a pentagonal ring of length nine,
contradicting $\tau(F)=10$. So $h'_3$ is a pentagon, and both $h'_1$
and $h'_4$ are hexagons. By Lemma \ref{lem4-3}(3), it follows that
$F$ is the graph (c) in Figure \ref{fig5-12}; that is also
$F^2_{36}$ in Figure \ref{fig5-20}.

{\em Case 2.} $\psi_{10}(F)=5$. By Eqs. (3) and (4), $n_5(R)=1$ and
$n_6(R)=5-\frac{1}{2}(5-r(R))$.

\begin{figure}[!hbtp]\refstepcounter{figure}\label{fig5-13}
\begin{center}
\includegraphics{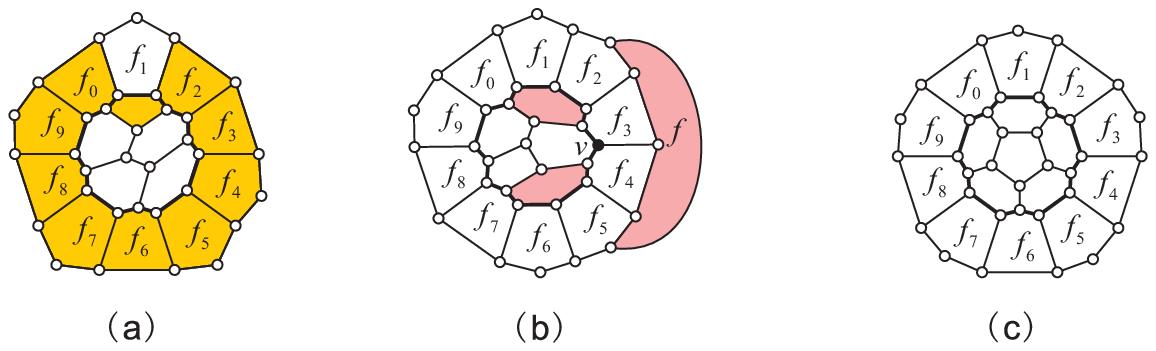}\\
{Figure \ref{fig5-13}: The illustration for Case 2 of the proof of
Lemma \ref{lem5-7}.}
\end{center}
\end{figure}

If there exist two vertices of $R$ on $C$ having a common neighbor
within $C$, let $G$ be the subgraph induced by $R$ together with all
vertices within $C$. By Lemma \ref{lem4-3}, it follows that $G$ is
isomorphic to the graph (a) or the graph (b) in Figure
\ref{fig5-13}. If $G$ is isomorphic to the graph (a), then $F$
contains a pentagonal ring $R'$ with length 10 and $s(R')=4$. Hence
$s(R')=4<\psi_{10}(F)\le s(R')$, a contradiction. So suppose $G$ is
isomorphic to the graph (b). Let $f$ be the face adjoining $R$ along
a 4-length path (see Figure \ref{fig5-13}, the common adjacent face
$f$ of $f_2,...,f_5$). If $f$ is a pentagon, then $F$ contains a
pentagonal ring with length 9 (the pentagonal ring
$(R-\{f_3,f_4\})\cup f$ in Figure \ref{fig5-13}), contradicting
$\tau(F)=10$. So suppose $f$ is a hexagon. Then $F$ has a  set
$\mathcal H$ of three mutually disjoint hexagons such that $f\in
\mathcal H$ and $F-\mathcal H$ has an isolated vertex $v$(the black
vertex in Figure \ref{fig5-13} (b)), contradicting that $F$ is
$3$-resonant. So there is no $k$-resonant fullerene graph satisfying
the condition.

Otherwise, any two 2-degree vertices on $C$ have distinct  neighbors
within $C$: If two 2-degree vertices on $C$ are adjacent, by Lemma
\ref{lem4-3}(2) the other three 2-degree vertices on $C$ would have
one common neighbor within $C$, a contradiction. Then the five faces
adjacent to $R$ within $C$ form a ring $R'$ of $F$ with $C$ as its
outer cycle. Note that the edges of $R'$ connecting  its outer cycle
and inner cycle form a cyclic 5-edge-cut of $F$. Since $|C|=15$, the
inner cycle of $R'$ bounds a face $f'$ of $F$ by Lemma \ref{lem5-6}.
So $f'$ is a pentagon. Therefore the subgraph $G$ induced by $R$
together with all vertices within $C$ is isomorphic to the graph (c)
in Figure \ref{fig5-13}. An analogous argument yields that the five
faces adjacent $G$ along its boundary are hexagons and $F-G$ is a
pentagon. So $F$ is $F_{40}$ as shown in Figure \ref{fig5-20}.

Finally  it suffices to prove that $F^2_{36}$ and $F_{40}$ are
$k$-resonant for every integer $k\geq 3$. For $F_{40}$,  it has a
perfect matching illustrated in Figure \ref{fig5-20} that is
alternating on  all hexagons. Hence $F_{40}$ is  $k$-resonant for
all $k\geq 1$. For $F^2_{36}$, let $G_1:=h_1\cup h_2\cup h_3\cup
h_4$ and $G_2:=h_1'\cup h_2'\cup h_4'\cup h_5'$. Then $G_1$ and
$G_2$ are disjoint, and the restrictions of  perfect matching $M$
illustrated in Figure \ref{fig5-20} on $G_1$ and $G_2$  are also
their perfect matchings. That means that the union of perfect
matchings of $G_1$ and $G_2$ can be extended to a perfect matching
of $F^2_{36}$. For each of $G_1$ and $G_2$, it is easy to see that
any disjoint hexagons are mutually resonant. Hence any disjoint
hexagons of $F^2_{36}$ forms a sextet pattern. \qed

\begin{figure}[!hbtp]\refstepcounter{figure}\label{fig5-15}
\begin{center}
\includegraphics{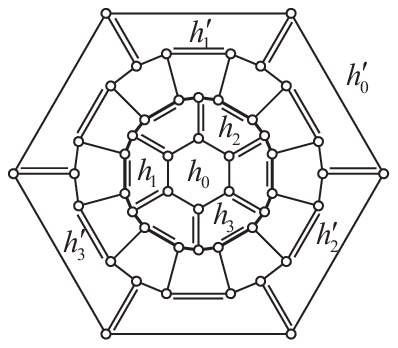}\\
{Figure \ref{fig5-15}: The  fullerene graph $F_{48}$  with a perfect
matching $M_0$.}
\end{center}
\end{figure}

\begin{lemma}\label{lem5-9}
Let $F$ be a  fullerene graph with $\tau(F)=12$. Then $F$ is
$3$-resonant if and only if $F$ is   $F_{48}$ shown in Figure
\ref{fig5-15}. Further, $F_{48}$ is $k$-resonant for any integer
$k\geq 3$.
\end{lemma}

\pf Let $F$ be a $3$-resonant fullerene graph with $\tau(F)=12$. Let
$R$ be the  pentagonal ring of $F$ with length $l(R)=12$ and $C$ the
inner cycle of $R$. Since $F$ has exactly 12 pentagons, there is no
pentagons within $C$ and hence $n_5(R)=0$. By Eqs (3) and (4),
$s(R)=6$ and $n_6(R)=4+\frac{1}{2}r(R)$. Let $v_0, v_1, \dots, v_5$
be the six 2-degree vertices of $R$  on $C$ arranged clockwise.

\begin{figure}[!hbtp]\refstepcounter{figure}\label{fig5-16}
\begin{center}
\includegraphics{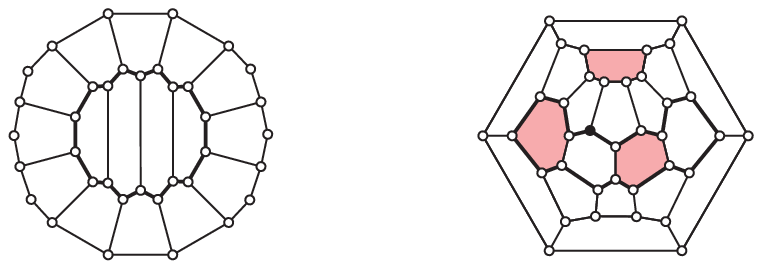}\\
{Figure \ref{fig5-16}: The subgraph $G_1$ (left) and the fullerene
graph containing $G_1$ (right).}
\end{center}
\end{figure}

If two 2-degree vertices on $C$ are connected by an edge of $F$
through the interior of $C$, then by Lemma \ref{lem4-3} (1) and (3),
the other four 2-degree vertices are connected by two edges of F
since every face within $C$ is a hexagon. Then  $r(R)=0$,  the
subgraph of $F$ induced by the vertices of $R$ is isomorphic to the
graph $G_1$ as shown in Figure \ref{fig5-16} (left) and $F$ is
 the fullerene graph as shown in Figure \ref{fig5-16}
(right). Then the  three  shadowed disjoint hexagons in Figure
\ref{fig5-16} (right)  are not mutually resonant, contradicting that
$F$ is $3$-resonant.

So we may suppose each 2-degree vertex on $C$ has a neighbor within
$C$. Then $r:=r(R)\geq 2$. Let $G$ be the subgraph induced by the
vertices within $C$ and  $m$  the number of  edges of $G$.
 Since $C$ has six 2-degree vertices, $3r-2m=6$.

If $G$ is not connected, then $G$ is a forest since $F$ is
cyclically 5-edge connected. Then $m=r-\omega$, where $\omega\ge 2$
is the number of components of $G$. So $3r-2(r-\omega)=6$, and $r=
6-2\omega\le 2$. Therefore $r=2$ and $G$ consists of two isolated
vertices, denoted by $u,v$. We may assume that
$N(u)=\{v_0,v_1,v_2\}$ and $N(v)=\{v_3,v_4,v_5\}$. Let $f$ be the
face within $C$ containing vertices $v_0,u,v_2,v_3,v$ and $v_5$.
Note that the six 2-degree vertices are not adjacent on $C$. So
$v_2v_3\notin E(f)$ and $v_5v_0\notin E(f)$. Therefore $|f|\ge 8$, a
contradiction.

So suppose that $G$ is connected. Let $\partial G$ be the boundary
of $G$ which is a closed walk. Note that a cut-edge of $G$ will
contribute 2 to $|\partial G|$. Let $x$ be the number of inner faces
of $G$. By Euler's formula, $m=r-1+x$. So $3r-2(r-1+x)=6$, and
 $r=4+2x$. On the other hand, since every inner face of $G$
is also a face of $F$, every inner face of $G$ is a hexagon. So
$|\partial G|+6x=2m=2r-2+2x=2(4+2x)-2+2x=6+6x$. Hence $|\partial
G|=6$.

\begin{figure}[!hbtp]\refstepcounter{figure}\label{fig5-17}
\begin{center}
\includegraphics{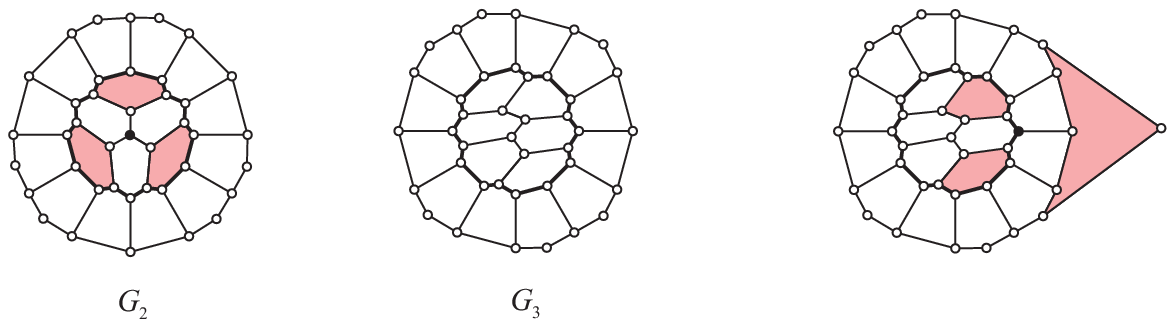}\\
{Figure \ref{fig5-17}: The subgraphs $G_2$  and  $G_3$.}
\end{center}
\end{figure}

If $x=0$, then $G$ is a tree. Since $|\partial G|=6$, $G$ has three
edges. So $G$ is isomorphic to a $K_{1,3}$ or a 3-length path. Hence
the subgraph of $F$ induced by $V(R\cup G)$ is isomorphic to either
$G_2$ or $G_3$ shown in Figure \ref{fig5-17}. Whenever $G_2\subset
F$ or $G_3\subset F$, $F$ has three  disjoint hexagons which are not
mutually resonant (see Figure \ref{fig5-17}).

Hence $x> 0$. Since the length of every cycle of $G$ is at least 6
and $|\partial G|=6$, $\partial G$ is a 6-length cycle. Since $F$ is
cubic, there are six edges connecting the 2-degree vertices
$v_0,v_1,\dots, v_5$ to vertices of $G$. So $G$ is a hexagon. By the
symmetry of $R$, $F$ is  the fullerene graph $F_{48}$ shown in
Figure \ref{fig5-15}.

Conversely, it suffices to show that $F_{48}$ is $k$-resonant for
$k\geq 3$. Let $M_0$ be the perfect matching of $F_{48}$ illustrated
in Figure \ref{fig5-15}. Let $f_0,f_1,f_2,f_3$ and
$f'_0,f'_1,f'_2,f'_3$ be the hexagons marked in $F_{48}$(see Figure
\ref{fig5-15}), and let $M_1:=M_0\oplus f_1\oplus f_2\oplus f_3$,
$M_2:=M_0\oplus f_1'\oplus f_2'\oplus f_3'$ and $M_3:=M_0\oplus
f_1\oplus f_2\oplus f_3\oplus f_1'\oplus f_2'\oplus f_3'$. Let
$\mathcal H$ be any set of mutually disjoint hexagons of $F_{48}$.
If $h_0,h'_0 \notin \mathcal H$, then every hexagon of $\mathcal H$
is $M_0$-alternating. If $h'_0 \notin \mathcal H$ but $h_0 \in
\mathcal H$, then  every hexagon of $\mathcal H$ is
$M_1$-alternating. If $h_0 \notin \mathcal H$ but $h'_0 \in \mathcal
H$, then  every hexagon of $\mathcal H$ is $M_2$-alternating. If
$\{h_0,h'_0\} \subseteq \mathcal H$, then $\mathcal H=\{h_0,h'_0\}$
and both $h_0$ and $h'_0$ are $M_3$-alternating. Thus $F_{48}$ is
$k$-resonant for $k\geq 3$. \qed

\vs

Summarizing the above results (Lemmas \ref{lem5-1}-\ref{lem5-9}), we
have the following main theorem.

\begin{thm}\label{thm3-10}
A  fullerene graph $F$ is $3$-resonant  if and only if $F$ is
 one of $F_{20}$, $F_{24}$, $F_{28}$, $F_{32}$,
$F_{36}^1$, $F_{36}^2$, $F_{40}$, $F_{48}$ and $\text{C}_{60}$.
Further, these nine fullerene graphs are all $k$-resonant for every
integer $k\geq 1$.
\end{thm}

From Theorem \ref{thm3-10}, we  arrive immediately at the the
following  result.

\begin{thm}\label{thm3-11}
A fullerene graph $F$ is $3$-resonant  if and only if it is
$k$-resonant for any integer $k\ge 3$.
\end{thm}

\section{Sextet polynomials of 3-resonant fullerene graphs}

The {\em sextet polynomial} of a benzenoid system $G$ for counting
sextet patterns was introduced by Hosoya and Yamaguchi \cite{HY} as
follows:
$$B_G(x)=\sum_{i=0}^{C(G)}\sigma(G,i)x^i,\eqno (7) $$
where $\sigma(G,i)$ denotes the number of sextet patterns of $G$
with $i$ hexagons, and $C(G)$ the Clar number of $G$. The sextet
polynomial of $\text{C}_{60}$ is computed \cite{SLZ1} as
$$B_{\text{C}_{60}}(x)=5x^8+320x^7+1240x^6+1912x^5+1510x^4+660x^3+160x^2+20x+1.\eqno (8)$$
For a detailed discussion and review of sextet polynomials, see
\cite{GU,R1}.

Since any independent hexagons of a 3-resonant fullerene graph form
a sextet pattern, we can compute easily the sextet polynomials of
the other eight $3$-resonant fullerene graphs as follows, by
counting sets of disjoint hexagonal faces.
$$\begin{array}{lll}
 B_{F_{20}}(x)&=&1,\\
B_{F_{24}}(x)&=&(x+1)^2=x^2+2x+1,\\
B_{F_{28}}(x)&=&(2x+1)^2=4x^2+4x+1,\\
B_{F_{32}}(x)&=&(3x+1)^2=9x^2+6x+1,\\
B_{F_{36}^1}(x)&=&2x^4+16x^3+20x^2+8x+1,\\
B_{F_{36}^2}(x)&=&(x^2+4x+1)^2=x^4+8x^3+18x^2+8x+1,\\
B_{F_{40}}(x)&=&(5x^2+5x+1)^2=25x^4+50x^3+35x^2+10x+1, \text{and}\\
B_{F_{48}}(x)&=&(2x^3+9x^2+7x+1)^2=4x^6+36x^5+109x^4+130x^3+67x^2+14x+1.
\end{array}
\eqno (9) $$

\section*{Acknowledgments}
The authors thank the anonymous referees for their helpful comments
on this paper and pointing out to us reference \cite{LKS}.


\end{document}